%% file: arxiv_aut.tex
\documentclass[final,5p,9pt,twocolumn,number,sort&compress]{elsarticle}{\tiny}
\input{preamble.tex}

\usepackage{enumitem}
\setitemize{itemsep=-1pt}

\journal{Automatica}

\begin{document}

\begin{frontmatter}


 \cortext[cor1]{This material is based upon work supported by the National Science Foundation under Grants 1728629 and 1728605. }

\title{Algorithms for Joint Sensor and Control Nodes Selection in Dynamic Networks}


\author[utsa]{Sebastian~A.~Nugroho}
\ead{sebastian.nugroho@my.utsa.edu}
\author[utsa]{Ahmad~F.~Taha}
\ead{ahmad.taha@utsa.edu}
\author[utsa]{Nikolaos~Gatsis}
\ead{nikolaos.gatsis@utsa.edu}
\author[utd]{Tyler~H.~Summers}
\ead{tyler.summers@utdallas.edu}
\author[utsa]{Ram~Krishnan}
\ead{ram.krishnan@utsa.edu}

\address[utsa]{Department of Electrical and Computer Engineering, The University of Texas at San Antonio, 1 UTSA Circle, San Antonio, TX  78249}
\address[utd]{Department of Mechanical Engineering, The University of Texas at Dallas, 800 W Campbell Rd, Richardson, TX, USA 75080}


\begin{abstract}
The problem of placing or selecting sensors and control nodes plays a pivotal role in the operation of dynamic networks. This paper proposes optimal algorithms and heuristics to solve the simultaneous sensor and actuator selection problem in linear dynamic networks. In particular, {\color{black} a sufficiency condition} of static output feedback stabilizability is used to obtain the minimal set of sensors and control nodes needed to stabilize an unstable network. We show the joint sensor/actuator selection and output feedback control can be written as a mixed-integer nonconvex problem.  
To solve this nonconvex combinatorial problem, three methods based on (1) mixed-integer nonlinear programming, (2) binary search algorithms, and (3) simple heuristics are proposed. The first method yields optimal solutions to the selection problem---given that some constants are appropriately selected. The second method requires a database of binary sensor/actuator combinations, returns optimal solutions, and necessitates no tuning parameters. The third approach is a heuristic that yields suboptimal solutions but is computationally attractive. The theoretical properties of these methods are discussed and numerical tests on dynamic networks showcase the trade-off between optimality and computational time. 
\end{abstract}


\end{frontmatter}


\section{Introduction}

\color{black}
\vspace{-0.2cm}
Consider an unstable dynamic network of $N$ interconnected nodes
	\begin{equation}~\label{equ:StateSpace-all-2}
		\m{\dot{x}}(t) = \m{Ax}(t)+\boldsymbol{Bu}(t), \;\;	\m{y}(t) = \m{Cx}(t)
	\end{equation}
 where $\m A$ has at least one unstable eigenvalue; $\m x(t), \m u(t),$ and $\m y(t)$ collect the state, input, and output vectors for all $N$ nodes. This paper studies the joint problems of \textit{(i)} stabilization of dynamic network \eqref{equ:StateSpace-all-2} through \textit{static output feedback control} (SOFC) while simultaneously \textit{(ii)} selecting or placing a minimal number of sensors and control nodes.

Problem \textit{(i)} corresponds to finding a control law $\m u(t)=\m F\m y(t)$ such that the closed loop system eigenvalues of $\m{A+BFC}$ are in the LHP~\citep{Astolfi2000}. This type of control is advantageous in the sense that it only requires output measurements rather than full state information, is analogous to the simple proportional controller, and can be implemented without needing an observer or an augmented dynamic system. Problem \textit{(ii)} corresponds to finding minimal number of \textit{sensors and actuators} (SA) yielding a feasible solution for the static output feedback (SOF) stabilization problem.  The joint formulations of Problems (\textit{i})--\textit{(ii)} can be abstracted through this high-level optimization routine: 
\begin{subequations}\label{eq:SAS-SOF}
\begin{align}
\min_{\m \Pi, \m \Gamma, \mF}  \;& \sum_{k=1}^{N} \pi_k + \gamma_k \\
	\st \; & \mathrm{real}(\mathrm{eig}(\mA+\mB \m \Pi \mF\m \Gamma \mC)) < 0,  \; \pi_i , \gamma_j  	\in \{0,1\}
\end{align}
\end{subequations}
where $\pi_i$ and $\gamma_j$ are binary variables selecting the $i$-th actuator and $j$-th sensor; $\m \Pi$ and $\m \Gamma$ are diagonal matrices containing all $\pi_i$ and $\gamma_j$. These binary variables post- and pre-multiply $\m B$ and $\m C$, thereby activating the optimal sensors and control nodes while designing a SOFC law. 
Even for small to mid-size dynamic networks, problem \eqref{eq:SAS-SOF} is difficult to solve as the SOFC problem---without the SA selection---is known to be nonconvex~\citep{Crusius1999} (presumed to be NP-hard \citep{PERETZ2016}), 
and the SA selection introduces binary variables thereby increasing the nonconvexity. To that end, the objective of this paper is to develop optimal algorithms and heuristics to solve Problem \eqref{eq:SAS-SOF}. Next, we summarize the recent literature on solving variants of~\eqref{eq:SAS-SOF}. 


Hundreds of studies have investigated the separate problem of minimally selecting/placing sensors \textit{or} actuators while performing state estimation \textit{or} state-feedback control. This paper, as mentioned above, studies the joint SA selection in the sense that an observer-based controller, which invokes the separation principle and requires a dynamic system module to perform state estimation, is not needed. For this reason, we do not delve into the literature of separate sensor or actuator selection.  Interested readers are referred to our recent work \citep{Taha2017d,Nugroho2018} for a summary on methods that solve the separate SA selection problems. The literature of addressing the \textit{simultaneous sensor and actuator selection problem} (SSASP) is summarized in what follows.


Several attempts have been made to address variants of the SSASP in dynamic networks through the more general \textit{dynamic output feedback control} (DOFC) framework. Specifically, the authors in~\citep{de2000linear} investigate the $\mathcal{H}_2$ minimization via DOFC with SA selection, in which a reformulated suboptimal problem in the form of mixed-integer semi-definite program (MI-SDP) is proposed and solved using a coordinate descent algorithm. 
In \citep{Argha2017}, the SSASP for multi-channel $\mc{H}_2$ DOFC with regional pole placement is addressed. In particular, the authors develop a semi-definite program (SDP) framework and propose a sparsity-promoting algorithm to obtain sparse row/column feedback control matrices. This approach ultimately yields binary SA selection, without needing binary variables. 
The same algorithm is then employed in \citep{Singh2018} for SSA selection with simpler $\mathcal{H}_2/\mathcal{H}_{\infty}$ formulations. The SA selection with control configuration selection problem is formulated in \citep{pequito2015} using structural design and graph theory, which is proven to be NP-hard. Although this particular problem is similar to the SSASP with SOFC given in \eqref{eq:SAS-SOF}, the problem proposed in \citep{pequito2015}, along with the algorithms, are based on the information of structural pattern of the dynamic matrix. The limitations of these studies are discussed next. 

\textcolor{black}{First, the majority of these works \citep{de2000linear,Argha2017,Singh2018} consider the $\mathcal{H}_2/\mathcal{H}_{\infty}$ control framework in conjunction with dynamic output feedback which requires an additional block of dynamical systems to construct the control action (which is not the case in SOFC). Second, the work in \citep{de2000linear} assumes that the number of SA to be selected is known \textit{a priori}, which for certain cases is not very intuitive. Third, the sparsity-promoting algorithm proposed in \citep{Argha2017,Singh2018} is based on convex relaxation of the $l_0$ norm---called re-weighted $l_1$ norm---which is then solved iteratively until the solution converges, thus making it not suitable for larger dynamic networks. The other drawback of this method is that arbitrary convex constraints on the binary selection variables are not easy to include. Finally, the algorithm proposed in \citep{pequito2015}---which interestingly runs in polynomial-time if the structure of the dynamic matrix is irreducible---only computes the structure and the corresponding costs of the feedback matrix (along with the sets of selected SA). }

As an alternative to the aforementioned methods, this paper proposes algorithms and heuristics to solve the SSASP for unstable dynamic networks via SOFC. Specifically, we use a sufficiency condition for SOFC from \citep{Crusius1999} which reduces the SOF control problem---without the SSA selection---from a nonconvex problem into a simple {linear matrix inequality} (LMI) feasibility problem. The developed approaches are based on MI-SDP, binary search algorithms, and simple heuristics that use the problem structure to find good suboptimal solutions. A preliminary version of this work appeared in~\citep{Nugroho2018} where we focus mainly on the MI-SDP approach. Here, we significantly extend this approach with the addition of binary search algorithms, heuristics, thorough analytical discussion of the properties of the developed methods, and comprehensive numerical experiments. The paper contributions and organization are discussed as follows. 
{\color{black}
	\begin{table}[t]
	\scriptsize		
	\caption{Acronyms used in the paper.}
	\label{tab:acro}
	\centering\begin{tabular}{ c|c }
		\hline 
		{DOFC} & Dynamic Output Feedback Control \\
	{LMI}/{BMI} & Linear/Bilinear Matrix Inequalities \\
	{MI-SDP} & Mixed Integer Semidefinite Program \\
				{SA} & Sensors and Actuators \\ 
	{SOFC} & Static Output Feedback Control \\
			{SSASP} & Simultaneous Sensor and Actuator Selection Problem \\ 
		\hline
	\end{tabular}
\vspace{-0.4cm}
\end{table}}

\noindent $\bullet$ Firstly, we formally introduce the SOF stabilizability problem (Section~\ref{sec:PbmForm}). The SSASP through SOFC is then formulated and shown to be a nonconvex problem with mixed-integer nonlinear matrix inequality (MI-NMI) constraints (Section~\ref{sec:PbmForm2}). We prove that the SSASP can be formulated as a MI-SDP, and the equivalence between the two is shown (Section~\ref{sec:Big-M}). The MI-SDP, if solved using combinatorial optimization techniques, yields an optimal solution to the SSASP. 

\noindent $\bullet$ As a departure from the MI-SDP approach, we introduce a routine akin to binary search algorithms that computes an optimal solution for SSASP---the proof of optimality is given. The routine requires a database of binary SA combinations (Section~\ref{sec:BSAlgo}). 

\noindent $\bullet$ A heuristic that scales better than the first two approaches is also introduced. The heuristic is based on constructing a simple logic of infeasible or suboptimal combinations of SA, while offering flexibility in terms of the tradeoff between the computational time and distance to optimality (Section~\ref{sec:HeuAlgo}). {\color{black} A brief discussion on the computational complexity as well as thorough numerical tests showcasing the applicability of the proposed algorithms are provided (Section~\ref{sec:complexity} and Section~\ref{sec:results}).}

The presented algorithms in this paper have their limitations which are all discussed with future work and concluding marks (Section \ref{sec:conc}).
\normalcolor

\vspace{-0.2cm}
\section{Notation}
\vspace{-0.2cm}

The symbols $\mathbb{R}^n$ and $\mathbb{R}^{p\times q}$ denote column vectors with $n$ elements and real-valued matrices with size $p$-by-$q$. The set of $n\times n$ symmetric and positive definite matrices are denoted $\mathbb{S}^{n}$ and $\mathbb{S}^{n}_{++}$. Italicized, boldface upper and lower case characters represent matrices and column vectors---$a$ is a scalar, $\m a$ is a vector, and $\m A$ is a matrix. Matrix $\m I_n$ is a $n\times n$ identity square matrix, while $\m 0$ and $\m O$ represent zero vectors and matrices of appropriate dimensions. For a square matrix $\m X$, the notation $\Lambda(\m X)$ denotes the set of all eigenvalues of $\m X$. The function $\mathrm{Re}(c)$ extracts the real part of a complex number $c$, whereas $\mathrm{Blkdiag}(\cdot)$ is used to construct a block diagonal matrix. For a matrix $\mX\in\mathbb{R}^{p\times q}$, the operator $\mathrm{Vec}(\mX)$ returns a stacked $pq\times 1$ column vector of entries of $\mX$, while $\mathrm{Diag}(\mY)$ returns a $n\times 1$ column vector of diagonal entries of square matrix $\mY\in\mathbb{R}^{n\times n}$. The symbol $\otimes$ denotes the Kronecker product. For any $x\in\mathbb{R}$, $\vert x\vert$ and $\lceil x \rceil$ denote the absolute value and ceiling function of $x$. The cardinality of a set $\mathcal{S}$ is denoted by $\vert\mathcal{S}\vert$, whereas $(0)^n$ denotes a $n$-tuple with zero-valued elements. 
\vspace{-0.1cm}
\section{Static Output Feedback Control Review}~\label{sec:PbmForm}
\vspace{-0.0cm}
In this section, we present some necessary background including the definition of static output feedback stabilizability given a fixed SA combination. This formulation is instrumental in deriving the SSASP.

Consider a dynamic network consisting of $N$ nodes/sub-systems with $\mathcal{N} = \lbrace 1,\ldots,N\rbrace$ defining the set of nodes. The network dynamics are given as:
\begin{subequations}~\label{equ:StateSpace-all}
	\begin{align}
		\m{\dot{x}}(t) &= \m{Ax}(t)+\boldsymbol{Bu}(t)~\label{equ:StateSpace-x}  \\
		\m{y}(t) &= \m{Cx}(t).~\label{equ:StateSpace-y} 
	\end{align}
\end{subequations}
The state, input, and output vectors corresponding to each node $i \in \mathcal{N}$ are represented by $\m{x}_i(t) \in \mathbb{R}^{n_{x_i}}$, $\m{u}_i(t) \in \mathbb{R}^{n_{u_i}}$, and $\m{y}_i(t) \in \mathbb{R}^{n_{y_i}}$. The global state, input, and output vectors are written as $\m{x}(t) \triangleq [\m{x}_1^{\top}(t),\ldots,\m{x}_N^{\top}(t)]^{\top}$, $\m{u}(t) \triangleq [\m{u}_1^{\top}(t),\ldots,\m{u}_N^{\top}(t)]^{\top}$, and $\m{y}(t) \triangleq [\m{y}_1^{\top}(t),\hdots,\m{y}_N^{\top}(t)]^{\top}$ where $\m{x}(t)\in \mathbb{R}^{n_{x}}$, $\m{u}(t)\in \mathbb{R}^{n_{u}}$, and $\m{y}(t)\in \mathbb{R}^{n_{y}}$. 
Without loss of generality, we assume that the input $\m{u}_i(t)$ and output $\m{y}_i(t)$ at each node $i$ only correspond to that particular node. The global input-to-state and state-to-output matrices can be constructed as $\m{B} \triangleq \mathrm{Blkdiag}(\m{B}_{1},\m{B}_{2},\ldots,\m{B}_{N})$ and $\m{C} \triangleq \mathrm{Blkdiag}(\m{C}_{1},\m{C}_{2},\ldots,\m{C}_{N})$ where $\m B \in \mathbb{R}^{n_x\times n_u} $ and $\m C \in \mathbb{R}^{n_y\times n_x} $. This assumption enforces the coupling among nodes to be represented in the state evolution matrix $\m{A} \in \mathbb{R}^{n_x\times n_x} $, which is realistic in various dynamic networks as control inputs and observations are often determined locally. Additionally, we also assume that $\m{B}$ and $\m{C}$ are full column rank and full row rank. This assumption eliminates the possibility of redundant control nodes and system measurements.  

The notion of SOF stabilizability for the dynamic system \eqref{equ:StateSpace-all} is provided first. Some needed assumptions are also given.

\begin{asmp}~\label{asmp:detectandstable}
	\color{black}
	The system \eqref{equ:StateSpace-all} satisfies the following conditions: (a) The pair $(\m A,\m B)$ is stabilizable; (b) the pair $(\m A, \m C)$ is detectable; (c) $\m{B}$ and $\m{C}$ are full rank.\normalcolor
\end{asmp}

\begin{mydef}~\label{def:OutputFeedbackStabilizability}
	The dynamical system \eqref{equ:StateSpace-all} is stabilizable via SOF if there exists $\m F \in \mathbb{R}^{n_u\times n_y}$ with control law given as $\m u(t) = \m F \m y(t)$ such that $\mathrm{Re}(\lambda) < 0$ for every $\lambda \in \Lambda(\m A_{\mathrm{cl}})$ where $\m A_{\mathrm{cl}} \triangleq \m{A}+\m{B}\m{F C}$.
\end{mydef}
The above definition and assumption are standard in the SOF control literature  \citep{syrmos1997static,KUCERA1995,rosinova2003necessary,Astolfi2000}. In this paper, we consider a simple yet well known necessary and sufficient condition for SOF stabilizability---given next.
\begin{myprs}[From \citep{syrmos1997static}]~\label{prs:OFS-BMI}
	The dynamical system \eqref{equ:StateSpace-all} is SOF stabilizable with output feedback gain $\m  F \in \mathbb{R}^{n_u\times n_y}$ if and only if there exists $\m P \in \mathbb{S}^{n_x}_{++}$ such that
	\begin{align}
		&\m{A}^{\top}\m P+\m{PA}+\m{C}^{\top}\m  F^{\top}\m B^{\top}\m P + \m P \m{B}\m{F C} 
		\prec 0.~\label{eq:OFS-BMI} 
	\end{align}
\end{myprs}

In Proposition \ref{prs:OFS-BMI}, the matrix inequality \eqref{eq:OFS-BMI} is nonconvex due to bilinearity in terms of $\m P$ and $\m  F$. 
{\color{black}Bilinear matrix inequalities (BMIs)} are of great interest to many researchers specifically in systems and control theory during the past decades \citep{Fukuda2001} because many control problems can be formulated as optimization problems with BMI constraints \citep{vanantwerp2000tutorial}. Indeed, many methods to address problems involving BMIs have been developed. Specifically, methods to solve the BMI for SOF stabilizability in a form similar to \eqref{eq:OFS-BMI} are proposed in \citep{Dinh2012,Hu2016}. These methods, based on successive convex approximation, linearize the BMI constraints around a certain strictly feasible point such that the nonconvex problem can be approximated by solving a sequence of convex optimization problems with {\color{black} LMI} constraints. Another less computationally intensive approach that can be used to address the SOF stabilizability problem is introduced in \citep{Crusius1999}. This approach allows the SOF problem to be solved in an LMI framework---presented in the following proposition.  
\begin{myprs}[From~\citep{Crusius1999}]~\label{prs:OFS}
	The dynamic network \eqref{equ:StateSpace-all} is SOF stabilizable if there exist $\m M \in \mathbb{R}^{n_u\times n_u}$, $\m P \in \mathbb{S}^{n_x}_{++}$, $\m N \in \mathbb{R}^{n_u\times n_y}$, and $\m  F \in \mathbb{R}^{n_u\times n_y}$ such that the following LMIs are feasible
	\begin{subequations}\label{eq:OFS}
		\begin{align}
			&\m{A}^{\top}\m P+\m{PA}+\m{C}^{\top}\m N^{\top}\m B^{\top}+\m{B}\m{N C} 
			\prec 0~\label{eq:OFS1} \\
			&\m B \m M = \m P \m B ~\label{eq:OFS2} 
		\end{align}
	\end{subequations}
	with SOF gain $\m{F}=\m{M}^{-1}\m{N}$. 
\end{myprs}
The conditions presented in Proposition~\ref{prs:OFS} are only sufficient for SOF stabilizability. Although this makes the SOF stabilization problem a lot easier to solve than the nonconvex BMI, this method has a drawback. In particular, the feasibility of solving \eqref{eq:OFS} as an SDP depends on the state-space representation of system \eqref{equ:StateSpace-all}, meaning that it is not guaranteed for \eqref{eq:OFS} to be feasible even if system \eqref{equ:StateSpace-all} is known to be stabilizable by SOF; see \citep{Crusius1999} for a relevant discussion. 
However, if the system is indeed stabilizable by SOF, then there exists a state-space transformation for system \eqref{equ:StateSpace-all} that leads \eqref{eq:OFS} to be feasible \citep{Crusius1999}. Further discussion about finding this transformation can be found in \citep{Crusius1999,neto1998}. With that in mind, and to develop tractable computational techniques to solve the SSASP, we use the LMI formulation of the SOF problem. The next section presents the problem formulation.
\section{Problem Formulation}~\label{sec:PbmForm2}
The  SSA selection through SOF control can be defined as the problem of jointly selecting a minimal set---or subset---of SA while still maintaining the stability of the system through SOF control scheme. To formalize the SSASP, let $\gamma_i\in \lbrace0,1\rbrace$ and $\pi_i \in \lbrace0,1\rbrace$ be two binary variables that represent the selection of SA at node $i$ of the dynamic networks. We consider that $\gamma_i = 1$ if the sensor of node $i$ is selected (or activated) and $\gamma_i = 0$ otherwise. Similarly, $\pi_i = 1$ if the actuator of node $i$ is selected and $\pi_i = 0$ otherwise. The augmented dynamics can be formulated as
\begin{subequations}\label{equ:StateSpaceSSA-all}
	\begin{align}
		\m{\dot{x}}(t) &= \m{Ax}(t)+\m{B\Pi u}(t) ~\label{equ:StateSpaceSSA-x}\\
		\m{y}(t) &= \m{\Gamma Cx}(t), ~\label{equ:StateSpaceSSA-y} 
	\end{align}
\end{subequations}
where $\m{\Pi}$ and $\m{\Gamma}$ are symmetric block matrices defined as 
\begin{subequations}\label{equ:defPiGamma}
	\begin{align}
		\m{\Pi} &\triangleq \mathrm{Blkdiag}(\pi_1\m{I}_{n_{u_{1}}},\pi_2\m{I}_{n_{u_{2}}},\ldots,\pi_N\m{I}_{n_{u_{N}}})\label{equ:defPi}\\
		\m{\Gamma}&\triangleq\mathrm{Blkdiag}(\gamma_1\m{I}_{n_{y_{1}}},\gamma_2\m{I}_{n_{y_{2}}},\ldots,\gamma_N\m{I}_{n_{y_{N}}}).\label{equ:defGamma}
	\end{align} 
\end{subequations}

The simultaneous sensor and actuator selection problem (SSASP) through SOF stabilizability can be written as in~\eqref{eq:OFS-SSAProblem}. 
The optimization variables of \eqref{eq:OFS-SSAProblem} are $\lbrace\m\pi,\m\gamma,\m N,\m M,\m P\rbrace$ with $\m\pi = [\pi_1,\pi_2,\ldots,\pi_N]^{\top}$, $\m\gamma = [\gamma_1,\gamma_2,\ldots,\gamma_N]^{\top}$; $\m{\Pi}$ and $\m{\Gamma}$ are the matrix forms of $\m\pi$ and $\m\gamma$ as defined in \eqref{equ:defPiGamma}. In the next sections, $\m\pi$ and $\m\gamma$ will be used interchangeably with $\m\Pi$ and $\m\Gamma$. Constraints \eqref{eq:OFS-SSAProblem-2} and \eqref{eq:OFS-SSAProblem-3} are obtained by simply applying the sufficient condition for SOF stabilizability. Constraint \eqref{eq:OFS-SSAProblem-5} is an additional linear logistic constraint which can be useful to model preferred  activation or deactivation of SA on particular nodes and to define the desired minimum and maximum number of activated SA. This constraint is also useful in multi-period selection problems where certain actuators and sensors are deactivated due to logistic constraints. The objective of the paper is to develop computational methods to solve SSASP. 

\begin{figure}[t]
	\vspace{-0.32cm}
	\centering	\includegraphics[scale=0.5]{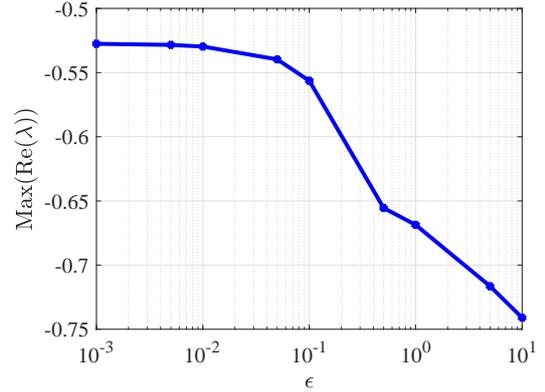}
	\vspace{-0.2cm}
	\caption{The relation between $\epsilon$ and maximum real closed loop system eigenvalues where $\lambda \in \Lambda(\m A+\mB\mF\mC)$.}
	\label{fig:0}
	\vspace{-0.38cm}
\end{figure}

\begin{mdframed}[style=MyFrame] 
	\vspace{-0.4cm} 
	\begin{subequations}~\label{eq:OFS-SSAProblem}
		\begin{align}
			\hspace{-0.0cm}\text{SSASP}\text{:}\;\; \minimize_{\substack{\m\pi,\m\gamma,\m N\\ \m M,\m P}} \;\; & \sum_{k=1}^{N}  
			\pi_k +  \gamma_k ~\label{eq:OFS-SSAProblem-1} \\
			\subjectto \;\; &\m{A}^{\top}\m P+\m{PA}+\m{C}^{\top}\m\Gamma\m N^{\top}\m\Pi\m B^{\top}\nonumber\\ &+\m{B}\m\Pi\m{N \Gamma C}\prec 0 ~\label{eq:OFS-SSAProblem-2}\\
			& \m{B}\m{\Pi M}=\m{P} \m B\m\Pi ~\label{eq:OFS-SSAProblem-3}\\\
			& \m{P} \succ 0~\label{eq:OFS-SSAProblem-4}\\
			&\m \Phi \bmat{\m \pi \\ \m \gamma} \leq \m \phi ~\label{eq:OFS-SSAProblem-5}\\
			&\m\pi \in \{0,1\}^N,\;\m\gamma \in \{0,1\}^N~\label{eq:OFS-SSAProblem-6}.
		\end{align}
	\end{subequations}
\end{mdframed}
\begin{rem}~\label{rem:marstab}
	The solution of \, {\normalfont \text{SSASP}} guarantees that the dynamic network is stabilized using the minimal number of SA, {\color{black} as the closed loop stability is ensured by the sufficient condition for the existence of SOFC given in Proposition \ref{prs:OFS}}. {\color{black} This entails that the closed-loop eigenvalues are in the left-half plane and close to the $j\omega$-axis. If it is desired to move the closed-loop eigenvalues further away from the $j\omega$-axis,} the matrix inequality~\eqref{eq:OFS-SSAProblem-2} can be upper bounded by $-\epsilon \m I$ where $\epsilon \in \mathbb{R}_{+}$. Larger values for $\epsilon$ will result in closed-loop eigenvalues further away from marginal stability. 
	Figure~\ref{fig:0} shows this relationship for a random dynamic network (described in Section~\ref{sec:results}  with $n_x=20, n_u=10,n_y=10$) given that all SA are activated after solving the LMI feasibility problem defined by \eqref{eq:OFS-SSAProblem-2} and \eqref{eq:OFS-SSAProblem-3} for different values of $\epsilon$. 
\end{rem}
{\color{black}
\begin{rem}
	It is important to mention that our focus here is to find a SOF control gain that stabilizes dynamical system \eqref{equ:StateSpaceSSA-all} with minimum number of SA. With that in mind, performance metrics such as robustness and energy cost functions are not considered in  {\normalfont \text{SSASP}}.
\end{rem}
}

After solving \eqref{eq:OFS-SSAProblem}, the selected SA are obtained and represented by $\lbrace\m \pi^*,\m \gamma^*\rbrace$. Due to SSA selection, the matrix $\mB\m\Pi$ will most likely not be full column rank, hence the existence of an invertible matrix $\mM$ is not assured. This is not the case when solving \eqref{eq:OFS} due to the fact that $\mB$ being full column rank  and $\mP \succ 0$ ensure $\mM$ to be nonsingular---see Lemma \ref{lem:0}. However, if \eqref{eq:OFS-SSAProblem} returns $\mM$ that is invertible, then the SSASP is solved with SOF gain $\mF$ to be computed as $\m{M}^{-1}\m{N}$. Otherwise, $\mF$ can be computed as $\hat{\m{M}}^{-1}\hat{\m{N}}$ where $\hat{\m{M}}$ and $\hat{\m{N}}$ are the submatrices of ${\m{M}}$ and ${\m{N}}$ that correspond to activated SA. Proposition \ref{prs:OFS-SOF} ensures the SOF stabilizbilty with minimal SA after solving SSASP. 

\begin{mylem}\label{lem:0}
Let $\m M \in \mathbb{R}^{m\times m}$ and $\m P \in \mathbb{S}^{n_x}$ be the solution of $\mB\mM = \mP\mB$ where $\m B \in \mathbb{R}^{n_x\times m}$ and $m \leq n_u$. If $\mP \succ 0$ and $\mathrm{Rank}(\mB) = m$, then $\mM$ is invertible. 
\end{mylem}

\begin{myprs}~\label{prs:OFS-SOF}
	Let $\m P$, $\m M$, $\m N$, $\m \Pi^*$, and $\m \Gamma^*$ be the solution of SSASP with appropriate dimensions. Also, let $\hat{\mB}\in \mathbb{R}^{n_x\times m}$, $\hat{\mC}\in \mathbb{R}^{r\times n_x}$, $\hat{\mM}\in \mathbb{R}^{m\times m}$, and $\hat{\mN}\in \mathbb{R}^{m\times r}$, where $m\leq n_u$ and $r\leq n_y$, be the matrices (or submatrices) representing the nonzero components of $\mB\m\Pi^*$, $\m\Gamma^*\mC$, $\m\Pi^*\mM$, and $\m\Pi^*\mN\m\Gamma^*$ that correspond to activated SA. Then, the closed loop system $\mA+\hat{\mB}\mF\hat{\mC}$ is stable with SOF gain $\mF = \hat{\mM}^{-1}\hat{\mN}$.
\end{myprs}

{\color{black} See \ref{apx:A} and \ref{apx:B} for the proofs of Lemma \ref{lem:0} and Proposition \ref{prs:OFS-SOF}}. SSASP \eqref{eq:OFS-SSAProblem} is nonconvex due to the presence of MI-NMI in the form of $\m\Pi\m{N \Gamma}$ and mixed-integer BMI in \eqref{eq:OFS-SSAProblem-3}. Therefore, it cannot be solved by any general-purpose mixed-integer convex programming solver. To that end, we present three approaches that solve or approximate \eqref{eq:OFS-SSAProblem}. The first approach---presented in the next section---is developed utilizing techniques from linear algebra and disjunctive programming principle~\citep{nemhauser1988integer,grossmann2002review}. The other two approaches are developed based on binary search algorithm and heuristics, as presented in the Sections~\ref{sec:BSAlgo} and~\ref{sec:HeuAlgo}.

\section{SSASP as a MI-SDP}~\label{sec:Big-M}
In this section, we present the first approach to solve~\eqref{eq:OFS-SSAProblem}, which transforms SSASP from a mixed-integer nonconvex problem to a MI-SDP. The following theorem presents this result. 
\color{black}
\begin{theorem}~\label{thrm:Big-M}
	SSASP is equivalent to
	\begin{mdframed}[style=MyFrame]
		\vspace{-0.37cm} 
		\begin{subequations}~\label{eq:Big-M-SSAProblem}
			\vspace{-0.1cm}\color{black}
			\begin{align}
				\minimize_{\substack{\m\pi,\m\gamma,\m N,\m M \\  \m P, \m \Theta}} & \;\;\;\sum_{k=1}^{N} \pi_k + \gamma_k ~\label{eq:Big-M-SSAProblem-1} \\
				\subjectto  & \notag\\
				& \hspace{-1.3cm}\m{A}^{\top}\m{P}+\m{PA}+\m{C}^{\top}\m\Theta^{\top}\m B^{\top}+\m{B}\m {\Theta C} \prec 0 ~\label{eq:Big-M-SSAProblem-2}\\
&\hspace{-1.3cm}\m\Psi_1(\mN,\m\Theta) \leq L_1\m\Delta_1(\m \Gamma,\m\Pi)~\label{eq:Big-M-SSAProblem-3} \\
&\hspace{-1.3cm}\m\Psi_2(\mM,\m \Omega(\mP)) \leq L_2\m\Delta_2(\m\Pi)~\label{eq:Big-M-SSAProblem-3-2}\\
&\hspace{-1.3cm}\m\Psi_3(\m \Xi(\mP)) \leq L_3\m\Delta_3(\m\Pi)~\label{eq:Big-M-SSAProblem-3-3}\\
				&\hspace{-1.73cm} 
			\;\;\;\;\;	\m{P} \succ 0, \;\; \m \Phi \bmat{\m \pi \\ \m \gamma} \leq \m \phi, \;\; \m\pi \in \{0,1\}^N,\;\m\gamma \in \{0,1\}^N,~\label{eq:Big-M-SSAProblem-4}
			\end{align}\normalcolor
		\end{subequations}\vspace{-0.3cm}
	\end{mdframed}
	where \eqref{eq:Big-M-SSAProblem-3},\eqref{eq:Big-M-SSAProblem-3-2},\eqref{eq:Big-M-SSAProblem-3-3} are linear constraints in which each function is specified as 
	{\small \begin{subequations} 
			\begin{align}
						\hspace{-1.4cm}	\m\Psi_1(\mN,\m\Theta) \triangleq &\bmat{\mathrm{Vec}(\m\Theta)\\-\mathrm{Vec}(\m\Theta)\\\mathrm{Vec}(\m\Theta)\\-\mathrm{Vec}(\m\Theta)\\\mathrm{Vec}(\m\Theta-\mN)\\-\mathrm{Vec}(\m\Theta-\mN)} \label{eq:bigM_aux_proof_1a}\\
				\m\Delta_1(\m \Gamma,\m\Pi) \triangleq  &\bmat{\mathrm{Diag}(\mI_{n_y}\otimes \m\Pi)\\\mathrm{Diag}(\mI_{n_y}\otimes \m\Pi)\\\mathrm{Diag}(\m\Gamma\otimes\mI_{n_u})\\\mathrm{Diag}(\m\Gamma\otimes\mI_{n_u})\\\mathrm{Diag}(2\mI_{n_u\times n_y}-\mI_{n_y}\otimes \m\Pi-\m\Gamma\otimes\mI_{n_u})\\\mathrm{Diag}(2\mI_{n_u\times n_y}-\mI_{n_y}\otimes \m\Pi-\m\Gamma\otimes\mI_{n_u})} \label{eq:bigM_aux_proof_1b}\\
				\m\Psi_2(\mM,\m \Omega(\mP))  \triangleq &\bmat{\mathrm{Vec}(\mM)\\-\mathrm{Vec}(\mM)\\\mathrm{Vec}(\m \Omega(\mP))\\-\mathrm{Vec}(\m \Omega(\mP))\\\mathrm{Vec}(\mM-\m \Omega(\mP))\\-\mathrm{Vec}(\mM-\m \Omega(\mP))} \label{eq:bigM_aux_proof_1c}\\
			\m\Delta_2(\m\Pi) \triangleq  &\bmat{\mathrm{Diag}(\mI_{n_u^2}-\mI_{n_u}\otimes \m\Pi+\m\Pi\otimes\mI_{n_u})\\\mathrm{Diag}(\mI_{n_u^2}-\mI_{n_u}\otimes \m\Pi+\m\Pi\otimes\mI_{n_u})\\\mathrm{Diag}(\mI_{n_u^2}+\mI_{n_u}\otimes \m\Pi-\m\Pi\otimes\mI_{n_u})\\\mathrm{Diag}(\mI_{n_u^2}+\mI_{n_u}\otimes \m\Pi-\m\Pi\otimes\mI_{n_u})\\\mathrm{Diag}(2\mI_{n_u^2}-\mI_{n_u}\otimes \m\Pi-\m\Pi\otimes\mI_{n_u})\\\mathrm{Diag}(2\mI_{n_u^2}-\mI_{n_u}\otimes \m\Pi-\m\Pi\otimes\mI_{n_u})}\label{eq:bigM_aux_proof_1d}\\
			\m\Psi_3(\m \Xi(\mP)) \triangleq &\bmat{\mathrm{Vec}(\m \Xi(\mP))\\-\mathrm{Vec}(\m \Xi(\mP))} \label{eq:bigM_aux_proof_1e}\\
			\m\Delta_3(\m\Pi) \triangleq  &\bmat{\mathrm{Diag}(\mI_{n_x\times n_u}-\m\Pi\otimes\mI_{n_x})\\\mathrm{Diag}(\mI_{n_x\times n_u}-\m\Pi\otimes\mI_{n_x})}\label{eq:bigM_aux_proof_1f},
			\end{align}
	\end{subequations}}
$\m \Omega$ and $\m \Xi$ are functions defined as
\begin{subequations}~\label{eq:Big-M-SSAProblem_aux_mat}
	\begin{align}
		\m \Omega(\mP) &\triangleq (\mB^{\top}\mB)^{-1}\mB^{\top}\mP\mB~\label{eq:Big-M-SSAProblem-new1}\\
		\m \Xi(\mP)&\triangleq (\mI-\mB(\mB^{\top}\mB)^{-1}\mB^{\top})\mP\mB,~\label{eq:Big-M-SSAProblem-new4}
	\end{align}	
\end{subequations}
	 $\m \Theta \in \mathbb{R}^{n_u\times n_y}$, is an additional optimization variable, and $L_1,L_2,L_3\in\mathbb{R}_{++}$ are predefined, sufficiently large constants. 
\end{theorem}
\normalcolor
\color{black}
\begin{rem}
	Although \eqref{eq:Big-M-SSAProblem} is equivalent to \,{\normalfont SSASP}, the quality of the solution that comes out of \eqref{eq:Big-M-SSAProblem} is very dependent on the choice of $L_{1,2,3}$. This observation is corroborated by numerical test results discussed in Section \ref{sec:results}.
\end{rem}
\normalcolor
The proof of Theorem \ref{thrm:Big-M} is given in \ref{apx:C}. This theorem allows the SSASP to be solved as a MI-SDP, which can be handled using a variety of optimization methods such as: branch-and-bound algorithms \citep{gamrath2016scip,gally2017}, outer approximations \citep{LubinYamangilBentVielma2016}, or cutting-plane methods \citep{Kong2010}. The next section presents a departure from MI-SDP to an algorithm that returns optimal solutions to SSASP, without requiring $L_{1,2,3}$.

\vspace{-0.0cm}
\section{Binary Search Algorithm for SSA Selection}\vspace{-0.3cm}~\label{sec:BSAlgo}

In this section we present an algorithm that is similar, in spirit, to binary search algorithms. The presented algorithm here seeks optimality for SSASP while not requiring any tuning parameters such as $L_{1,2,3}$; see Theorem~\ref{thrm:Big-M}.  

\vspace{-0.3cm}
\subsection{Definitions and Preliminaries}
In what follows, we provide some needed definitions.
\color{black}
\begin{mydef}~\label{def:setS}
	Let $\mathcal{S}_\pi$ and $\mathcal{S}_\gamma$ be two $N$-tuples representing the selection of SA, i.e., $\mathcal{S}_\pi \triangleq (\pi_1,\ldots,\pi_N)$ and $\mathcal{S}_\gamma \triangleq (\gamma_1,\ldots,\gamma_N)$. Then, the selection of SA can be defined as $\mathcal{S} \triangleq (\mathcal{S}_\pi,\mathcal{S}_\gamma)$ such that $\lbrace\m \Pi,\m \Gamma\rbrace = \mathcal{G}(\mathcal{S})$, $\m\Pi = \mathcal{G}_{\pi}(\mathcal{S})$, and $\m\Gamma = \mathcal{G}_{\gamma}(\mathcal{S})$ where $\mathcal{G}:\{0,1\}^{2N}\rightarrow \mathbb{R}^{n_u\times n_u}\times \mathbb{R}^{n_y\times n_y}$, $\mathcal{G}_{\pi}:\{0,1\}^{2N}\rightarrow \mathbb{R}^{n_u\times n_u}$, and $\mathcal{G}_{\gamma}:\{0,1\}^{2N}\rightarrow \mathbb{R}^{n_y\times n_y}$ are linear mappings. The number of nodes with activated SA is defined as $\mathcal{H}(\mathcal{S}) \triangleq \sum_{k=1}^{N} \pi_k + \gamma_k$ where $\mathcal{H}:\mathcal{S}\rightarrow \mathbb{Z}_{+}$.
\end{mydef}
\normalcolor
\color{black}
\begin{mydef}~\label{def:bigS}
	Let $\mathbfcal{S}\triangleq \lbrace\mathcal{S}_q\rbrace_{q=1}^\sigma$ be the candidate set such that it contains all possible combinations of SA where $\sigma$ denotes the number of total combinations, i.e., $\sigma\triangleq\vert\mathbfcal{S}\vert$. Then, the following conditions hold: 
	\begin{enumerate}
		\item For every $\mathcal{S}\in\mathbfcal{S}$, it holds that 
		\begin{align}
		\mathcal{S}\in\left\{ \mathcal{S}\in \{0,1\}^{2N} \,|\, \mathcal{G}(\mathcal{S}) \textit{ is feasible for \eqref{eq:OFS-SSAProblem-5}}\right\}, \nonumber
		\end{align}
		\item For every $q$ where $1\leq q \leq \sigma$, $\mathcal{H}(\mathcal{S}_{q-1}) \leq \mathcal{H}(\mathcal{S}_q)$.
	\end{enumerate}
\end{mydef}
\normalcolor
\begin{mydef}~\label{def:tripletS}
	For any $\mathcal{S}_q\in\mathbfcal{S}$ such that $\lbrace\m \Pi_q,\m \Gamma_q\rbrace = \mathcal{G}(\mathcal{S}_q)$, $\m B_q$ and $\m C_q$ are defined as the matrices containing the nonzero components of $\m B \m \Pi_q$ and $\m \Gamma_q\m C$ that correspond to the activated SA. Then, we say that $\mathcal{S}_q$ is feasible for \eqref{eq:OFS} if and only if the triplet $(\mA,\m B_q,\m C_q)$ is feasible for \eqref{eq:OFS}.
\end{mydef}

The following example shows how the candidate set $\mathbfcal{S}$ is constructed for a given simple logistic constraint.
\color{black} 
\vspace{-0.1cm}
\setlength{\textfloatsep}{8pt}
{\small \begin{algorithm}[h]
	\caption{\small Binary Search Algorithm (BSA)}\label{alg:DaC}
	\DontPrintSemicolon
	\textbf{initialize:} $\mathcal{S}^*=(1)^{2N}$, $p=1$\;
	\textbf{input:}  $\mathbfcal{S}_p$, $\mA$, $\mB$, $\mC$\;
	\While{$\mathbfcal{S}_p \neq \emptyset$}{
		\textbf{compute:} $\sigma\leftarrow\vert\mathbfcal{S}_p\vert$, $q \leftarrow \lceil \sigma/2 \rceil$, $\mathcal{S}_q\in\mathbfcal{S}_p$\;
		\eIf{$\mathcal{S}_q$ is feasible for \eqref{eq:OFS}\label{step:step5}}{
			$\mathcal{S}^*\leftarrow\mathcal{S}_q$, $\mathbfcal{S}_p \leftarrow \mathbfcal{S}_p\setminus\lbrace\mathcal{S}\in\mathbfcal{S}_p\,\vert\,\mathcal{H}(\mathcal{S}) \geq \mathcal{H}(\mathcal{S}_q)\rbrace$\;
		}{
			$\mathbfcal{S}_p\leftarrow\mathbfcal{S}_p\setminus\lbrace\mathcal{S}\in\mathbfcal{S}_p\,\vert\,\mathcal{S}_q\vee\mathcal{S} = \mathcal{S}_q \rbrace$ \label{alg1:step9}\;
		}
		$p \leftarrow p + 1$\;
	}
	\textbf{output:} $\mathcal{S}^*$\;
\end{algorithm} }\setlength{\floatsep}{8pt}
\normalcolor
\vspace{-0.1cm}
\color{black}
\begin{exmpl}\label{ex:Ex1}
	Suppose that the dynamical system has two nodes. If the logistic constraint dictates that $ 1 \leq \mathcal{H}(\mathcal{S}) = \sum_{k=1}^{2} \pi_k + \gamma_k <  4$ for all $\mathcal{S}\in\mathbfcal{S}$, then $\mathbfcal{S}$ can be constructed as
	\begin{align*}
		\mathbfcal{S} = \Big\lbrace
		&(1,0,0,0),(0,1,0,0),(0,0,1,0),(0,0,0,1),\\&(1,1,0,0),(1,0,1,0),(1,0,0,1),(0,1,1,0),\\&(0,1,0,1),(0,0,1,1),(1,1,1,0),(1,1,0,1),\\&(1,0,1,1),(0,1,1,1)\Big\rbrace.
	\end{align*}
\end{exmpl} 
\normalcolor
\vspace{-0.6cm}
\subsection{Binary Search Algorithm to Solve SSASP}\vspace{-0.0cm}\label{BSA-1}

The objective of this algorithm is to find \textit{an} optimal solution $\mathcal{S}^*\in\mathbfcal{S}$ such that $\mathcal{H}(\mathcal{S}^*)\leq\mathcal{H}(\mathcal{S})$ for all $\mathcal{S}\in\mathbfcal{V}$ where $\mathbfcal{V}\triangleq\lbrace \mathcal{S}\in\mathbfcal{S} \,|\, \mathcal{S}\textit{ is feasible for \eqref{eq:OFS}}\rbrace$---{\color{black}that is, for any $\mathcal{S}\in\mathbfcal{V}$, there exist a corresponding feedback gain $\mF$ that stabilizes dynamical system \eqref{equ:StateSpaceSSA-all}}. Realize that any  $\mathcal{S}_q\in\mathbfcal{V}$, $(\mA,\m B_q,\m C_q)$ is feasible for SSASP with objective function value equal to $\mathcal{H}(\mathcal{S}_q)$. 

The routine to solve SSASP based on binary search algorithm is now explained and summarized in Algorithm~\ref{alg:DaC}. Let $p$ be the index of iteration and $q$ be the index of  position in the ordered set $\mathbfcal{S}$. Hence at iteration $p$, the candidate set containing all possible combinations of SA can be represented as $ \mathbfcal{S}_p$, with $\sigma=\vert\mathbfcal{S}_p\vert$, and any element of $ \mathbfcal{S}_p$ at position $q$ can be represented by $\mathcal{S}_q$. Also, let $\mathcal{S}^*$ be the known best solution at iteration $p$. 

Next, obtain $\mathcal{S}_q$ where $\mathcal{S}_q\in\mathbfcal{S}_p$ and $q = \lceil\sigma/2\rceil$. At this step, we need to determine whether system \eqref{equ:StateSpaceSSA-all} is SOF stabilizable with the particular combination of SA $\lbrace\m \Pi_q,\,\m \Gamma_q\rbrace = \mathcal{G}(\mathcal{S}_q)$. To that end, we solve the LMIs~\eqref{eq:OFS}. If $\mathcal{S}_q$ is feasible for \eqref{eq:OFS}, then update $\mathcal{S}^*\leftarrow\mathcal{S}_q$. Since $\mathcal{S}_q$ is feasible, then we can discard all combinations that have more or equal number of activated SA, i.e., the combinations that are suboptimal. Otherwise, if $\mathcal{S}_q$ is infeasible for \eqref{eq:OFS}, $\mathcal{S}_q$ can be discarded along with all combinations that (a) have less number of activated SA than $\mathcal{S}_q$ \textbf{and} (b) the activated SA are included in $\mathcal{S}_q$. Realize that the above method reduces the size of $\mathbfcal{S}_p$ in every iteration since one or more elements of $\mathbfcal{S}_p$ are discarded. The algorithm now continues and terminates whenever $\mathbfcal{S}_p = \emptyset$. The details of this algorithm are given in Algorithm~\ref{alg:DaC}. Example \ref{exp:2} gives an illustration how $\mathbfcal{S}_p$ is constructed in every iteration.

\begin{exmpl}\label{exp:2}
	Consider again the  dynamic system from Example~\ref{ex:Ex1}. Let $(1,0,0,1)$ be the starting combination and, for the sake of illustration, assume that \eqref{eq:OFS} is infeasible for this combination. Then, by Algorithm \ref{alg:DaC}, combinations $(1,0,0,0)$ and $(0,0,0,1)$ are discarded. The candidate set now comprises the following elements 
	\begin{align*}
		\mathbfcal{S}_2= \Big\lbrace
		&(0,1,0,0),(0,0,1,0),
		(1,1,0,0),(1,0,1,0),\\&(0,1,1,0),(0,1,0,1),(0,0,1,1),(1,1,1,0),\\&(1,1,0,1),(1,0,1,1),(0,1,1,1)\Big\rbrace.
	\end{align*}
	Let $(0,1,0,1)$ be the new starting point and assume that this combination is feasible for \eqref{eq:OFS}. Then, all combinations that have greater or equal number of activated SA can be discarded. The remaining possible candidates in the candidate set are
	\begin{align*}
		\mathbfcal{S}_3 = \Big\lbrace
		&(0,1,0,0),(0,0,1,0)\Big\rbrace.
	\end{align*}
	This algorithm continues in a fashion similar to the above routine. If none of these combinations in $\mathbfcal{S}_3$ is feasible, then Algorithm~\ref{alg:DaC} returns $\mc{S}^*= (0,1,0,1)$ as the solution.  
\end{exmpl} 

In what follows, we discuss the optimality of Algorithm~\ref{alg:DaC} through Theorem~\ref{thm:thm2}---see \ref{apx:D} for the proof.

\begin{theorem}\label{thm:thm2}
	Algorithm~\ref{alg:DaC} yields an optimal solution of SSASP.
\end{theorem}

\begin{rem}
	$\mathcal{S}^*$ from Algorithm~\ref{alg:DaC} might not be unique. This is the case since there could be multiple binary combinations of SA yielding the same number of activated SA and hence the same objective function value {\color{black}$\sum_{k=1}^{N} \pi_k + \gamma_k$.}
\end{rem}
\vspace{-0.1cm}

\section{Heuristics to Solve SSASP}~\label{sec:HeuAlgo}
%
%
The binary search algorithm in the previous section requires the construction of the candidate set $\mathbfcal{S}$ in an off-line database, while leading to an optimal solution for the SSASP. Seeking optimality and constructing an off-line database might be impractical for large-scale dynamic networks. 
Moreover, the other approach presented in Section~\ref{sec:Big-M} entails solving \eqref{eq:Big-M-SSAProblem}, a MI-SDP, which might consume large computational resources. This motivates the development of a heuristic for the SSASP that forgoes optimality---the focus of this section.

In short, the heuristic builds a dynamic, virtual database of all possible combinations---not by generating all of these combinations, but by having a procedure that identifies suboptimal/infeasible candidates---while attempting to find a SA combination that has the least number of activated SA that makes system \eqref{equ:StateSpaceSSA-all} SOF stabilizable. The high-level description of this heuristic is given as follows:
\begin{enumerate}
	\item Generate a random SA string $\mathcal{S}$;
	\item If $\mathcal{S}$ is in the \textit{forbidden set} (the set of suboptimal/infeasible combinations), repeat Step 1;
	\item If $\mathcal{S}$ is infeasible for~\eqref{eq:OFS}, add $\mathcal{S}$ to {\color{black}the forbidden set} and repeat Step 1;
	\item If $\mathcal{S}$ is feasible for~\eqref{eq:OFS}, add $\mathcal{S}$ to a set of candidate strings and discard suboptimal candidates, repeat Step 1 with fewer activated SA;
	\item Three metrics guide how many times these steps are repeated: $\mathrm{maxRandom}$, $\mathrm{maxInfeasibility}$, $\mathrm{maxIter}$.
\end{enumerate}

{ \begin{algorithm}[t]
	\caption{\small Heuristic (HEU) to Solve SSASP}\label{alg:Heuristic}
	\DontPrintSemicolon
	\textbf{initialize:}  $\mathcal{S}^* = (1)^{2N}$, $\mathbfcal{Z} = \mathbfcal{W}$, $\barbelow{w}$, and $\bar{w}$\;
	\textbf{set:} $t = 1$, $p = 1$, $q = \lceil (\barbelow{w}+\bar{w})/2\rceil$\;
	\textbf{input:} $\mathrm{maxIter}$, $\mathrm{maxInfeasibility}$\;
	\While{$p \leq \mathrm{maxIter}$ {\bf and} $\barbelow{w} \leq q \leq \bar{w}$ }
	{
		\While{$p \leq \mathrm{maxIter}$ {\bf and} $t \leq \mathrm{maxInfeasibility}$}
		{
			\textbf{compute:} $\mathcal{S}^{(q)}_p \notin \mathbfcal{Z}$ from Algorithm~\ref{alg:RandomCandidate}\label{step:step6}\;
			\eIf {$\mathcal{S}^{(q)}_p\neq(0)^{2N}$ \label{step:step8}}
			{
				\eIf {\eqref{eq:OFS} is feasible}
				{
					$\mathcal{S}^*\leftarrow\mathcal{S}^{(q)}_p $, $\bar{w}\leftarrow q-1$, $t\leftarrow 1$, $p\leftarrow p+1$\;
					\textbf{break}\;
				}{
					$\mathbfcal{Z}\leftarrow\mathbfcal{Z}\,\cup\,\lbrace\mathcal{S}^{(q)}_p\rbrace$, $t\leftarrow t+1$, $p\leftarrow p+1$\;
				}
			}{
				$t\leftarrow 1$\;
				\textbf{break}\;
			}
		}
		\eIf {$t > \mathrm{maxInfeasibility}$}
		{
			$q \leftarrow \lceil (q+\bar{w})/2\rceil$, $t\leftarrow 1$\;
		}{
			$q \leftarrow \lceil (\barbelow{w}+\bar{w})/2\rceil$\;
		}
	}
	\textbf{output:} $\mathcal{S}^*$\;	
\end{algorithm}}

We now introduce the details of the heuristic. Define \color{black} $$\mathbfcal{W}\triangleq \left\{\mathcal{S}\in \{0,1\}^{2N} \,|\, \mathcal{G}(\mathcal{S}) \textit{ is not feasible for \eqref{eq:OFS-SSAProblem-5}}\right\},$$that is, $\mathbfcal{W}$ is a finite set that comprises all combinations of SA that do not satisfy the logistic constraint \eqref{eq:OFS-SSAProblem-5}. Since we are interested in finding a candidate $\mathcal{S}$ that is feasible for \eqref{eq:OFS}, all elements in $\mathbfcal{W}$ do not need to be known. Instead, we just need to check whether $\mathcal{S}\notin\mathbfcal{W}$ using the logistic constraint \eqref{eq:OFS-SSAProblem-5}. Next, from the logistic constraint \eqref{eq:OFS-SSAProblem-5}, we define $\barbelow{w}\triangleq \mathrm{min}(\mathcal{H}(\mathcal{S}))$ and $\bar{w}\triangleq \mathrm{max}(\mathcal{H}(\mathcal{S}))$ for all $\mathcal{S}\notin\mathbfcal{W}$. That is, $\barbelow{w}$ and $\bar{w}$ represent the required minimum and maximum number of activated SA so that any candidate $\mathcal{S}\notin\mathbfcal{W}$ must satisfy $\barbelow{w}\leq \mathcal{H}(\mathcal{S}) \leq\bar{w}$. More importantly, $\barbelow{w}$ and $\bar{w}$ can also be used to bound the search space of a potential candidate $\mathcal{S}$. In contrast with Algorithm \ref{alg:DaC}, the heuristic constructs and updates---in each iteration---a set that contains combinations of SA that are known to be infeasible for \eqref{eq:OFS-SSAProblem}. This finite set is referred to as the \textit{forbidden set} and symbolized by $\mathbfcal{Z}$. Clearly, $\mathbfcal{W}\subseteq\mathbfcal{Z}$. Thus, any candidate $\mathcal{S}$ must not belong in $\mathbfcal{Z}$ because any $\mathcal{S}\in\mathbfcal{Z}$  is infeasible for \eqref{eq:OFS} and/or $\mathcal{S}\in\mathbfcal{W}$. To get a potential candidate $\mathcal{S}$ for this heuristic, we can randomly generate $\mathcal{S}$ such that $\mathcal{S}\notin\mathbfcal{Z}$.  

{ \begin{algorithm}[t]
	\caption{\small Candidate Generation}\label{alg:RandomCandidate}
	\DontPrintSemicolon
	\textbf{initialize:} $\mathcal{S}^{(q)}_p=(0)^{2N}$, $r = 1$\;
	\textbf{input:} $p, q, \barbelow{w}$\;
	\While{$r\leq \mathrm{maxRandom}$}{
		Randomly generate $\mathcal{S}$ with $\mathcal{H}({\mathcal{S}}) = q$\;
		\eIf{$\mathcal{S}\notin \mathbfcal{Z}$}{
			$\mathcal{S}^{(q)}_p \leftarrow \mathcal{S}$\;
			\textbf{break}\;
		}{
			$r \leftarrow r + 1$\;
			\If{$r > \mathrm{maxRandom}$}{
				$\barbelow{w} \leftarrow q + 1$\;
			}
		}
	}
	\textbf{output:} $\mathcal{S}^{(q)}_p, \barbelow{w}$\;
\end{algorithm}}

The heuristic is described as follows and summarized in Algorithm \ref{alg:Heuristic}. An essential part of the algorithm is a simple procedure to obtain a candidate $\mathcal{S}$ such that $\mathcal{S} \notin \mathbfcal{Z}$---this procedure is shown in Algorithm \ref{alg:RandomCandidate}. First, from the logistic constraint, $\barbelow{w}$, $\bar{w}$, and $\mathbfcal{Z}$ are initialized. Let $p$ denote the iteration index and $q = \lceil(\barbelow{w}+\bar{w})/2\rceil$ denote the desired number of activated SA for the candidate $\mathcal{S}$ such that $\mathcal{H}({\mathcal{S}}) = q$. Then, a candidate at iteration $p$ with $q$ number of activated SA can be denoted by $\mathcal{S}^{(q)}_p$. The next step is to generate a candidate $\mathcal{S}^{(q)}_p$ such that $\mathcal{S}^{(q)}_p \notin\mathbfcal{Z}$. As mentioned earlier, one simple method to obtain $\mathcal{S}^{(q)}_p$ is to randomly generate $\mathcal{S}$ with $q$ number of activated SA such that $\mathcal{S}\notin\mathbfcal{Z}${\color{black}---see Algorithm \ref{alg:RandomCandidate} for the detailed steps.} If such candidate cannot be obtained after some combinations of SA have been randomly generated, then $\mathcal{S}^{(q)}_p \leftarrow(0)^{2N}$. When this happens, we can assert that the majority of combinations of SA with $q$ or less than $q$ number of activated SA most likely belong to the forbidden set $\mathbfcal{Z}$. Given this condition, the required minimum number of activated SA can then be increased and updated. 
{\color{black}
If $\mathcal{S}^{(q)}_p$ is nonzero, then we must check whether $\mathcal{S}^{(q)}_p$ is feasible for \eqref{eq:OFS}. Then, if $\mathcal{S}^{(q)}_p$ is feasible for \eqref{eq:OFS}, we update $\mathcal{S}^*\leftarrow\mathcal{S}^{(q)}_p$; otherwise, update the forbidden set so that $\mathbfcal{Z}\leftarrow\mathbfcal{Z}\,\cup\,\lbrace\mathcal{S}^{(q)}_p\rbrace$.} Unlike Algorithm \ref{alg:DaC}, here we define $\mathrm{maxInfeasibility}$ that allows \eqref{eq:OFS} to be solved repeatedly with different candidates while having the same number of activated SA. This process is repeated until there exists a candidate that makes \eqref{eq:OFS} feasible or $\mathrm{maxInfeasibility}$ is reached. If \eqref{eq:OFS} is still infeasible, we increase the required number of activated SA for the next candidate, hoping that adding more activated SA will increase the chance for \eqref{eq:OFS} being feasible. The algorithm continues and terminates when maximum iteration, denoted by $\mathrm{maxIter}$, is reached or there is no more candidates that can be generated. At the end of Algorithm \ref{alg:Heuristic}, the best suboptimal combination of SA is given as $\mathcal{S}^*$. 

The algorithm in its nature allows the trade-off between the computational time and distance to optimality. This trade-off can be designed via selecting large values for $\mathrm{maxRandom},$ $\mathrm{maxInfeasibility},$ and $\mathrm{maxIter}$. The parameter $\mathrm{maxIter}$ depends on how the user is willing to wait before the algorithm terminates, $\mathrm{maxRandom}$ imposes an upper bound on how many times a random SA candidate $\mathcal{S}$ is generated such that it is does not belong to the forbidden set. Finally, $\mathrm{maxInfeasibility}$ defines how many LMI feasibility problems are solved with a fixed number of activated SA. 

\begin{table*}[t]
	\scriptsize
	\centering\color{black}
	\begin{threeparttable}[b]
	\caption{  Numerical test results on random network with $N = 10$, $n_x=20, n_u=10$ and  $n_y=20$. The symbol $\Delta t(s)$ represents the computational time measured in seconds. All methods successfully return a stable closed loop system, albeit close to the $j\omega$-axis; see Remark~\ref{rem:marstab}.The number of iterations for the MI-SDP corresponds to the BnB solver in YALMIP. 
	}
	\renewcommand{\arraystretch}{1.5}
	\begin{tabular}{|c|c|c|c|c|c|}
		\hline
		Scenario & {$\mathrm{Max}(\mathrm{Re}(\Lambda(\m A + \m B\m \Pi^*\mF \m \Gamma^*\m C)))$} & {$\Delta t(s)$} & {$\mathrm{Iterations}$} & {$\sum_k \pi_k + \gamma_k$} & $\m{\gamma}^*$ and $\m{\pi}^*$ \\
		\hline \hline 
		\multirow{2}[1]{*}{MI-SDP-1} & \multirow{2}[1]{*}{$-4.78 \times 10^{-3}$} & \multirow{2}[1]{*}{$24.2$} & \multirow{2}[1]{*}{$187$} & \multirow{2}[1]{*}{$3$} & $\m{\gamma}^*= \{0,     0,     0,     0,     0,     0,     1,     0,     0,     0\}$ \\
		&       &       &       &       & $\m{\pi}^*= \{   0,     0,     0,     0,     0,     1,     0,     0,     0,     1\}$ \\
		\hline
		\multirow{2}[1]{*}{MI-SDP-2} & \multirow{2}[1]{*}{$-6.71\times 10^{-3}$} & \multirow{2}[1]{*}{$13.3$} & \multirow{2}[1]{*}{$82$} & \multirow{2}[1]{*}{$5$} & $\m{\gamma}^* = \{1,  1,     0,     0,     0,     0,     0,     0,     0,     0\}$ \\
		&       &       &       &       & $\m{\pi}^*= \{   0,     1,     0,     0,     1,     0,     0,     0,     0,     1\}$ \\
		\hline
		\multirow{2}[0]{*}{BSA} & \multirow{2}[0]{*}{$-2.69\times 10^{-3}$} & \multirow{2}[0]{*}{$61.8$} & \multirow{2}[0]{*}{$198$} & \multirow{2}[0]{*}{$3$} & $\m{\gamma}^*= \{ 0,	0,	0,	0,	1,	0,	0,	0,	0,	0\}$ \\
		&       &       &       &       & $\m{\pi}^*= \{ 1,	0,	0,	0,	0,	1,	0,	0,	0,	0\}$ \\          \hline
		\multirow{2}[0]{*}{HEU\tnote{$\dagger$}} & \multirow{2}[0]{*}{$-1.65\times 10^{-3}$} & \multirow{2}[0]{*}{{$27.4$}} & \multirow{2}[0]{*}{{$50$}} & \multirow{2}[0]{*}{{$3.42$\tnote{$\ddagger$}}} & $\m{\gamma}^*= \{1,	0,	0,	0,	0,	0,	0,	0,	0,	0 \}$ \\
		&       &       &       &       & $\m{\pi}^*= \{1,	0,	0,	0,	0,	1,	0,	0,	0,	0\}$ \\
		\hline
	\end{tabular}%
	\label{tab:table1}%
	   \begin{tablenotes}
		\item[$\dagger$] The displayed values are mean values of 500 randomizations. The corresponding binary configurations of SA are taken from the first randomization. All of the 500 randomizations return stable closed loop system.
		\item[$\ddagger$] Out of the 500 randomizations, 294 return 3 activated SA, 198 return 4, and 8 randomizations return 5.
	\end{tablenotes}\normalcolor
	\end{threeparttable}
		\vspace{-0.1cm}
\end{table*}%

\vspace{-0.2cm}
\section{Discussion on Computational Complexity}\label{sec:complexity}
\vspace{-0.1cm}

\color{black}

{To discuss the computational complexity of the developed approaches, we start by discussing that of SDPs/LMIs. Primal-dual interior-point methods for SDPs have a worst-case complexity estimate of $\mathcal{O}\left(m^{2.75}L^{1.5} \right)$, where $m$ is the number of variables
	and $L$ is the number of constraints \citep{boyd1994linear}. The number of variables and constraints in~\eqref{eq:OFS} are}
\begin{subequations}
	\begin{align}
	m &=\underbrace{0.5 n_x(n_x+1)}_{\text{entries in}\,\,\m P} + \underbrace{n_u n_y}_{\text{entries in}\,\,\m N} + \underbrace{n_u^2}_{\text{entries in}\,\,\m M} \\
	L &=\underbrace{0.5 n_x(n_x+1)}_{\text{matrix inequality constraints}} +  \underbrace{n_x n_u}_{\text{ equality constraints}}.
	\end{align}
\end{subequations}

In various problems arising in control systems studies, it is shown that the complexity estimate is closer to $\mathcal{O}\left(m^{2.1}L^{1.2} \right)$ which is significantly smaller than the worst-case estimate $\mathcal{O}\left(m^{2.75}L^{1.5} \right)$~\citep{boyd1994linear}. Hence, the complexity estimate for the LMIs is $\mathcal{O}\left(n_x^{4.2}n_x^{2.4} \right)=\mathcal{O}\left(n_x^{6.6} \right)$, since typically $n_x > n_u$ and $n_x > n_y$ in dynamic networks. Admittedly, these figures are old now as many advancements in interior-point methods are often implemented within newer versions of SDP solvers. Given that, the first approach in Theorem~\ref{thrm:Big-M} entails solving the MI-SDP. 
Unfortunately, the worst-case complexity of solving MI-SDPs through branch-and-bound is $\mathcal{O}\left(2^{2N} n_x^{6.6} \right)$ as there are $2N$ SSA decision variables (for a network comprising $N$ nodes). However, as branch-and-bound solvers almost always terminate way before trying all combinations, it is very difficult to obtain the best-case performance for this approach. 

{The second approach, namely Algorithm~\ref{alg:DaC}, entails solving a LMI feasibility problem at each iteration. The best case performance of this algorithm occurs when a feasible solution is always obtained for the LMI at each iteration of Algorithm~\ref{alg:DaC}. This yields a logarithmic reduction in the number of candidate optimal solutions. Hence, the best case complexity of running Algorithm~\ref{alg:DaC} is of $\mathcal{O}\left( \log (2^{2N})\, n_x^{6.6} \right)=\mathcal{O}\left(N\, n_x^{6.6} \right)=\mathcal{O}\left(n_x^{7.6} \right)$. The worst-case performance of Algorithm~\ref{alg:DaC} occurs when the LMI returns an infeasible solution at each iteration, which is hard to quantify. Unfortunately, it is virtually impossible to upper or lower bound the number of these combinations as this is system-dependent.  Algorithm~\ref{alg:Heuristic} also entails either solving an LMI feasibility problem, while depending on a maximum iteration number and thresholds. The computational complexity hence depends on the user-defined maximum iteration number and thresholds which is hard to estimate here.}
\normalcolor

\section{Numerical Experiments}\vspace{-0.1cm}~\label{sec:results}
\color{black}
Numerical experiments are presented here to tests the proposed approaches on two dynamic networks. The first system is a random dynamic network adopted from \citep{Motee2008,MihailoSiteUnstableNodes}, whereas the second is a mass-spring system \citep{lin2013design}. Both systems are initially unstable and the latter has a sparser structure than the former. 
All simulations are performed using \iftrue MATLAB R2016b running on a 64-bit Windows 10 with 3.4GHz Intel Core i7-6700 CPU and 16 GB of RAM. All optimization problems are solved using MOSEK version 8.1~\citep{mosekAps} with YALMIP \citep{lofberg2004yalmip}.  All the MATLAB codes used in the paper are available for the interested reader upon request.
\normalcolor

\color{black}
\subsection{Comparing the Proposed Algorithms}
In the first part of our numerical experiment, we focus on testing the performance of the proposed approaches to solve SSASP on a relatively small dynamic network where optimality for the SSASP can be determined. Specifically, we consider the aforementioned random dynamic network with 10 subsystems, with two states per subsystem, so that 10 sensors and 10 actuators are available ($n_x=20, n_u = 10, n_y = 20$). Each sensor measures the two states per subsystem. We impose a logistic constraint so that there are at least one sensor and one actuator to be activated: $\sum_{i}^N\pi_i \geq 1$ and $\sum_j^N \gamma_j \geq 1$. In this particular experiment, we consider the following scenarios.
\normalcolor
\begin{enumerate}
	\item \textit{MI-SDP-1}:
	The first scenario uses the results from Theorem~\ref{thrm:Big-M} that shows the equivalence between SSASP and \eqref{eq:Big-M-SSAProblem}---the latter is solved via YALMIP's MI-SDP branch and bound (BnB) solver~\citep{yalmipMI}. We choose $L_1= 10^4$, $L_2 = 5\times 10^{6}$, and $L_3 = 5\times 10^{6}$. 
	The maximum number of iterations of the BnB solver is chosen to be 1000.
	\item \textit{MI-SDP-2}: The second scenario is identical to the first one with the exception that $L_1 = 10^4$, $L_2 =  10^{7}$, and $L_3 =  10^{7}$. This scenario shows the impact of tuning parameters $L_{1,2,3}$ on the performance of the MI-SDP approach.
	\item \textit{BSA}: The third scenario directly follows Algorithm \ref{alg:DaC} and solves \eqref{eq:OFS} in each iteration to check the feasibility of the given SA combinations, while also computing the SOF gain matrix simultaneously from the solution of LMI~\eqref{eq:OFS}. 
	\item \textit{HEU}: In the fourth scenario, we implement the heuristic as described in Algorithms \ref{alg:Heuristic} and \ref{alg:RandomCandidate}. The parameters of the heuristic in this scenario are $\mathrm{maxRandom} = 10^4$, $\mathrm{maxInfeasibility} = 10$, and $\mathrm{maxIter} = 50$. Since the proposed heuristic entails randomizations to generate SA candidates, we perform 500 randomizations. Then, from these randomizations, the mean values are computed.
\end{enumerate}

Table \ref{tab:table1} presents the result of this test. All scenarios return solutions with stable closed-loop system eigenvalues---albeit close to marginal stability (see Remark~\ref{rem:marstab}).  In this simulation we are only concerned with strict stabilization through SOF control---the focus here is the minimal number of activated SA and computational time.
The MI-SDP approach yields the optimal solution as discussed in Theorem~\ref{thrm:Big-M} and confirmed by BSA and Theorem~\ref{thm:thm2}. Specifically, Theorem~\ref{thm:thm2} shows that the BSA returns the optimal solution to SSASP, and Theorem~\ref{thrm:Big-M} shows the equivalence between SSASP---a mixed integer problem with nonlinear matrix inequalities---and~\eqref{eq:Big-M-SSAProblem}, a MI-SDP. While MI-SDP-1 returns an optimal solution with a smaller computational time in comparison with BSA, the former is dependent on the choice of $L_{1,2,3}$. MI-SDP-2 solves the same problem for different values of $L_{1,2,3}$ and yields 5 activated SA---clearly a suboptimal solution. 
BSA here is advantageous in the sense that it does not require tuning to find the appropriate constants $L_{1,2,3}$. Also, BSA requires only an LMI solver, instead of a BnB solver for the MI-SDP. 

Out of the 500 randomizations for the HEU, 294 randomizations return 3 activated SA, 198 return 4, and 8 randomizations return 5 activated SA---yielding an average of 3.42 activated SA as shown in Table~\ref{tab:table1}. The HEU hence surprisingly yields optimal solutions with 3 activated SA in the majority of randomizations, while requiring a smaller computational time in comparison with BSA 
and a comparable computational time with the MI-SDP while not requiring any tuning for $L_{1,2,3}$. This implies that the BSA (optimal) and the heuristic (suboptimal) can be used for general dynamic networks with minimal testing, i.e., given the triplet $(\m A,\m B, \m C)$ these algorithms can be immediately applied. In contrast, the MI-SDP approach can be used if optimal solutions for the SSASP are sought when constructing an offline database for the BSA is impractical due to limitations on storage capacity. These contrasts are illustrated in the next numerical experiment.
\begin{table}[t]
\scriptsize	
	\centering\color{black}
	\begin{threeparttable}[b]
	\caption{Results for the MI-SDP and heuristic with $N = 50$, $n_x=100, n_u=50, n_y=100$. }
	\renewcommand{\arraystretch}{1.4}
\begin{tabular}{|c|c|c|c|}
	\hline
	{Method}&{$\mathrm{Max}(\mathrm{Re}(\Lambda(\m A + \m B\m \Pi^*\mF \m \Gamma^*\m C))$} & {$\Delta t(s)$} & {$\sum_{k} \pi_k + \gamma_k$}\\ \hline\hline
	MI-SDP & $-7.59\times 10^{-3}$& $13749.9$ & $9$ \\ \hline
	Heuristic &  ${-1.68\times 10^{-2}}$ & $6181.5$ & ${17.39}$\tnote{$\S$} \\ \hline
\end{tabular}~\label{table2}
\vspace{-0.2cm}
\begin{tablenotes}
	\item[$\S$] Out of the 10 randomizations, 4 return 3 activated SA while the remaining 6 returns various number of activated SA each ranging from 11 to 23.
\end{tablenotes}\normalcolor
\end{threeparttable}
		\vspace{-0.2cm}
\end{table}%


\color{black}
Secondly, to find out how our methods perform on a larger system, we test our MI-SDP and heuristic (Algorithm~\ref{alg:Heuristic}) to solve a larger random dynamic network consisting 50 subsystems with $n_x=100$, $n_u = 50$, and $n_y = 100$.  For MI-SDP, we set a maximum of 1000 iterations for the BnB solver.
As for heuristic, we perform 10 randomizations with $\mathrm{maxRandom} = 10^6, \mathrm{maxInfeasibility} = 10$, and $\mathrm{maxIter} = 200$. 
The result of this experiment is summarized in Table \ref{table2}.
 The MI-SDP via BnB algorithm returns 9 activated SA in total (3 sensors and 6 actuators) while requiring 112 iterations and $1.38\times10^{4}\, \mathrm{s} \approx 3.8 \,\mathrm{hrs}$ of computational time which we presume is the optimal solution. On the other hand, we obtain an average of $17.39$ activated SA with $6.18\times10^{3}\, \mathrm{s} \approx 1.7 \,\mathrm{hrs}$ of average computational time. Realize that the total available SA are 50 sensors and 50 actuators, which implies that there are $2^{100} \approx 1.26 \times 10^{30}$ combinations of SA. With the average number of activated SA being relatively small compared to the number of available SA, we conclude that the heuristic produces a reasonably good solution---in comparison with the MI-SDP solution.
 This experiment demonstrates the tradeoff between MI-SDP and heuristic. 
\normalcolor

\color{black}
\subsection{Comparative Study with Dynamic Output Feedback}
In the second part of the numerical experiment, we consider the mass-spring systems from \citep{lin2013design} of various sizes to measure the performance of our heuristic (Algorithm~\ref{alg:Heuristic}) relative to the other method by comparing it with the sparsity promoting algorithm (SPA) for $\mathcal{H}_{\infty}$ dynamic output feedback developed in \citep{Singh2018}. The heuristic algorithm is configured in a way such that $\mathrm{maxRandom} = 10^6$, $\mathrm{maxInfeasibility} = 10$, and $\mathrm{maxIter} = 200$. For each number of nodes $N$, we perform 10 randomizations for the heuristic in which the mean value of maximum real part of closed loop eigenvalues, computational time, and the number of activated SA are computed accordingly. The SPA is set up so that the maximum iteration number is 50 and the convergence tolerance is 0.5. The results are given in Table \ref{tab:mass_spring_heuristics_spa}.

\begin{table}[t]
	\vspace{-0.1cm}
	\centering\color{black}
	\begin{threeparttable}[b]
	\renewcommand{\arraystretch}{1.5}
	\caption{Numerical comparison results between the heuristic and SPA  with mass-spring systems. The function $\bar{\lambda}_{\mathrm{Re}}(\cdot)$ computes the maximum value of the real part of eigenvalues of $\m A_{cl}$ or $\tilde{\m A}_{cl}$, where $\m A_{cl}=\m A + \m B\m \Pi^*\mF \m \Gamma^*\m C$ and $\tilde{\m A}_{cl}$ represents the overall closed-loop dynamics matrix of the plant with DOFC. The notation $\m \Sigma$ represents the number of activated SA, i.e., $\sum_{k=1}^{N} \pi_k + \gamma_k$.}
	\vspace{0.2cm}
	\centering	\scriptsize	
	\begin{tabular}{|c|c|c|c|c|c|c|}
\hline
\multirow{2}{*}{$\m N$} & \multicolumn{3}{c|}{Heuristic} & \multicolumn{3}{c|}{SPA} \\ \cline{2-7} 
                   & $\bar{\lambda}_{\mathrm{Re}}(\m A_{cl})$      & $\Delta {t} (s)$   & $\m \Sigma$ & $\bar{\lambda}_{\mathrm{Re}}(\tilde{\m A}_{cl})$      &  $\Delta {t} (s)$      & $\m \Sigma$     \\ \hline\hline
$10 $    &    $-3.3\times 10^{-3}$   &  $6.9$      &  $2$     &       $-2.7\times 10^{-1}$ &    $10.4$    &   $10$    \\ \hline
$20 $  &    $-4.6\times 10^{-4}$    &  $19.5$      &  $2$     &       $-3.4\times 10^{-1}$ &     $129.1$   &   $20$    \\ \hline
$30$&    $-1.5\times 10^{-4}$    &  $98.9$      &   $2$    &       $-8.2\times 10^{-2}$ &   $1383.4$     &   $24$    \\ \hline
  $40 $ &  $-6.1\times 10^{-5}$      &  $406.0$      & $2$     &       $-2.6\times 10^{-1}$ &   $3844.7$     &  $40$     \\ \hline
 $50 $&  $-2.3\times 10^{-5}$      &  $1736.7$      &  $2.30$     &       $-5.8\times 10^{-2}$ &   $28513.9$     &   $26$    \\ \hline
\end{tabular}\label{tab:mass_spring_heuristics_spa}
\end{threeparttable}
\normalcolor
\vspace{0.1cm}
\end{table}

From this experiment, we observe that both methods are able to give stable closed-loop systems with the heuristic returning fewer activated SA than SPA. The computational time of the heuristic is significantly faster than SPA's. For example, if we consider the particular case of $N = 50$, the heuristic returns $2.30$ of activated SA on average, while SPA returns 26. 
These findings, however, should not conflated with the objectives of the $\mathcal{H}_{\infty}$ and SOFC methods, seeing that both methods consider different control metrics (static output feedback considers pure stabilization whereas dynamic output feedback with $\mathcal{H}_{\infty}$ control considers robustness) with different problem size and complexity. The results shown here are meant to  give the indication that when robustness is considered as a metric through dynamic output feedback control and SA selection, the corresponding problem requires more computational time, the activation of more SA, yet returns a closed-loop system that is more robust to disturbances. 

\normalcolor

\color{black}
\section{Summary, Limitations, and Future Directions}\vspace{-0.3cm}~\label{sec:conc}

In this paper we propose computational methods to solve the simultaneous SA selection problem (SSASP) through static output feedback control framework. Three different approaches to obtain the minimal selection of activated SA that yield stable closed-loop systems are proposed. The first approach utilizes disjunctive programming principles and linear algebra techniques to convert the mixed-integer nonconvex problem into a MI-SDP. The second approach uses a simple algorithm that is akin to the binary search algorithm. The third approach is a simple heuristic that constructs a dynamic data structure with infeasible combinations. 

The first two approaches yield optimal solutions for the {SSASP}, while the third yields suboptimal results while resulting in an improved computational time. In particular, the first approach requires finding suitable constants, namely $L_{1,2,3}$, to obtain an optimal solution to the MI-SDP and hence the {SSASP}. This is for combinatorial optimization problems that use the Big-M method or the McCormick Relaxation. The second optimal approach requires efficient data structures to store and update the feasible combinations of sensors and actuators, without requiring any tuning parameters. The third approach is suitable for larger networks as it trades optimality with improved computational time. It is noteworthy to mention that the second and third algorithms only require an LMI solver making it easier to interface with without the need to install any additional optimization packages or tune any parameters. 

The limitations of the proposed methods in this paper are listed as follows. 
%
First, we do not consider any robustness or energy metrics through SOFC and SSASP. For example, an interesting extension can capture the minimal SA selection alongside designing energy-aware and robust output feedback control laws for the activated nodes. Second, our approach still requires solving LMIs. Albeit the LMIs solved in this paper are simple with few optimization variables and only one block, this approach does \textit{not} scale graciously when large-scale dynamic networks with tens of thousands of nodes are studied.  

To that end, our future work in this topic will focus on addressing the aforementioend paper's limitations by deriving SOFC energy and robustness metrics with the SSASP, and consequently examining the performance of the presented algorithms in this paper. Furthermore, and to address the computational burden of solving many LMIs, we plan to investigate algebraic conditions on the existence of SOFC given a fixed SA selection. The idea here is to avoid solving LMI feasibility problem at each iteration. Instead, we can learn the feasibility of specific SA selection using these algebraic conditions. Finally, and instead of solving the MI-SDP form of {SSASP} using the classical BnB algorithm, we plan to test the performance of the outer approximations \citep{LubinYamangilBentVielma2016} and the cutting-plane methods \citep{Kong2010}---as these approaches have shown significant savings for the computational time especially when compared with classical BnB methods.



\normalcolor

\bibliographystyle{abbrv} 
\bibliography{bibliography}

\appendix
\section{Proof of Lemma \ref{lem:0}}\label{apx:A}
\color{black}
\begin{proof}
	Let $\m M$, $\m P$, and $\m B$ satisfy $\mB\mM = \mP\mB$. Then, suppose that $\mP \succ 0$ and $\mathrm{Rank}(\mB) = m$. Since $\mB$ is full column rank, $\mB^{\top}\mB$ is nonsingular. Premultiplying both sides of $\mB\mM = \mP\mB$ with $\mB^{\top}$ yields
	\begin{align*}
		\mB^{\top}\mB\mM = \mB^{\top}\mP\mB\;\Rightarrow\;\mM = (\mB^{\top}\mB)^{-1}\mB^{\top}\mP\mB.
	\end{align*}
	Since $\mP\succ 0$ and $\mathrm{Rank}(\mB) = m$, by \citep[Observation 7.1.8]{horn2012matrix}, we have $\mB^{\top}\mP\mB\succ 0$. Thus $\mM^{-1} = (\mB^{\top}\mP\mB)^{-1}\mB^{\top}\mB$.
\end{proof}
\normalcolor
\section{Proof of proposition \ref{prs:OFS-SOF}}\label{apx:B}
\color{black}
\begin{proof}
	Without loss of generality, $\mB\m\Pi^*$ and $\mC^{\top}\m\Gamma^*$ can be expressed as $\mB\m\Pi^* = \bmat{\hat{\mB}&\mO}$ and  $\mC^{\top}\m\Gamma^* = \bmat{\hat{\mC}^{\top}&\mO}$. Then, $\mM$ and $\mN$ can be partitioned as
	\begin{align*}
		\mM = \bmat{\mM_1&\mM_2\\\mM_3&\mM_4},\quad \mN = \bmat{\mN_1&\mN_2\\\mN_3&\mN_4},
	\end{align*}
	where $\mM_1\in \mathbb{R}^{m\times m}$ and $\mN_1\in \mathbb{R}^{m\times r}$. By letting $\hat{\mN} = \mN_1$, \eqref{eq:OFS-SSAProblem-2} can be expressed as
	\begin{align}
		\m{A}^{\top}\m P+\m{PA}+\bmat{\hat{\mC}^{\top}&\mO}\bmat{\hat{\mN}^{\top}&\mN_3^{\top}\\\mN_2^{\top}&\mN_4^{\top}}\bmat{\hat{\mB}^{\top}\\\mO} &\nonumber  \\
		+\bmat{\hat{\mB}&\mO}\bmat{\hat{\mN}&\mN_2\\\mN_3&\mN_4}\bmat{\hat{\mC}\\\mO}&\prec 0 \nonumber \\ 
		\Leftrightarrow\m{A}^{\top}\m P+\m{PA}+\hat{\mC}^{\top}\hat{\mN}^{\top}\hat{\mB}^{\top}+\hat{\mB}\hat{\mN}\hat{\mC}&\prec 0.~\label{eq:prs4-1}
	\end{align}
	Since \eqref{eq:OFS-SSAProblem-2} is feasible for $\m P$ and $\mN$, then \eqref{eq:prs4-1} is also feasible. Similarly, by letting $\hat{\mM} = \mM_1$, \eqref{eq:OFS-SSAProblem-3} can be expressed as
	\begin{align}
		\bmat{\hat{\mB}&\mO}\bmat{\hat{\mM}&\mM_2\\\mM_3&\mM_4} &= \mP\bmat{\hat{\mB}&\mO} \nonumber \\
		\Leftrightarrow \bmat{\hat{\mB}\hat{\mM}&\hat{\mB}\mM_2} &= \bmat{\mP\hat{\mB}&\mO}.~\label{eq:prs4-2}
	\end{align}
	Realize that \eqref{eq:prs4-2} holds since we assume that \eqref{eq:OFS-SSAProblem-3} holds. From \eqref{eq:prs4-2}, we have $\hat{\mB}\hat{\mM}=\mP\hat{\mB}$. Since $\hat{\mB}$ is full column rank and $\mP \succ 0$, by Lemma \ref{lem:0}, $\hat{\mM}$ is nonsingular. Finally, having $\m P\succ 0$, $\hat{\mM}$, and $\hat{\mN}$ that satisfy \eqref{eq:prs4-1} and $\hat{\mB}\hat{\mM}=\mP\hat{\mB}$, then according to Proposition \ref{prs:OFS}, the closed loop system $\mA+\hat{\mB}\mF\hat{\mC}$ is stable with $\mF = \hat{\mM}^{-1}\hat{\mN}$.
\end{proof}
\normalcolor
\section{Proof of Theorem \ref{thrm:Big-M}}\label{apx:C}
\color{black}
\begin{proof}
	Let $(\m{\Pi N \Gamma})_{ij}$ be the $(i,j)$ element of $\m{\Pi N \Gamma}$ and $\left\lbrace\pi_i,\,\gamma_j\right\rbrace$ be the associated SA selection that corresponds to $(\m{\Pi N \Gamma})_{ij}$. Then, there exists $\m \Theta \in \mathbb{R}^{n_u\times n_y}$ such that $\m \Theta = \m{\Pi N \Gamma}$. This relation is established as follows. Realize that, as $\m{\Pi}$ and $\m{\Gamma}$ are symmetric diagonal matrices with binary values, we can write $(\m{\Pi N \Gamma})_{ij}$ as
	\[
	(\m{\Pi N \Gamma})_{ij}=
	\begin{cases}
	\m{N}_{ij}, &\text{if}\;\;\pi_i\wedge\gamma_j = 1 \\
	0, & \text{if}\;\;\pi_i\wedge\gamma_j = 0.
	\end{cases}
	\]
	for $i = 1,\hdots,n_u$ and $j = 1,\hdots,n_y$. That is, if  $\left\lbrace\pi_i,\,\gamma_j\right\rbrace = \left\lbrace1,\,1\right\rbrace$, then $\m N_{ij} = \m \Theta_{ij}$. Otherwise, $\m N_{ij} \in\mathbb{R}$ and $\m \Theta_{ij} = 0$. For an appropriate large constant $L_1$, this is equivalent to 
	\begin{align*}
		\vert\m\Theta_{ij}\vert &\leq L_1 \pi_i,\;
		\vert\m\Theta_{ij}\vert \leq L_1 \gamma_j,\;
		\vert\m\Theta_{ij}-\m N_{ij}\vert \leq L_1(2-\pi_i-\gamma_j),
	\end{align*}
	which can be represented as $\m\Psi_1(\mN,\m\Theta) \leq L_1\m\Delta_1(\m \Gamma,\m\Pi)$ where $\m\Psi_1(\mN,\m\Theta)$ and $\m\Delta_1(\m \Gamma,\m\Pi)$ are given in \eqref{eq:bigM_aux_proof_1a} and \eqref{eq:bigM_aux_proof_1b}.
	This establishes \eqref{eq:Big-M-SSAProblem-3}. Consequently, $\m{\Pi N \Gamma}$ in \eqref{eq:OFS-SSAProblem-2} can be replaced with $\m{\Theta}$, giving \eqref{eq:Big-M-SSAProblem-2}. Now, pre-multiplying both sides of \eqref{eq:OFS-SSAProblem-3} with $\mB\mB^{\dagger}$ and, since $\mB\mB^{\dagger}\mB = \mB$, we obtain  
	\begin{subequations}
		\begin{align}
			\mB\mB^{\dagger}\mB\m\Pi\mM &= \mB\mB^{\dagger}\mP\mB\m\Pi \label{proof1_0} \\
			\Leftrightarrow \mB\m\Pi\mM &= \mB\mB^{\dagger}\mP\mB\m\Pi.\label{proof1_1a}
		\end{align}	
	\end{subequations}
	Since $\mB$ is full column rank, $\mB^{\dagger}$ can be expressed as $\mB^{\dagger}=(\mB^{\top}\mB)^{-1}\mB^{\top}$. Then by the same reason, \eqref{proof1_1a} implies 
	\begin{align}
		\m\Pi\mM = \mB^{\dagger}\mP\mB\m\Pi 
		\Leftrightarrow \m\Pi\mM = (\mB^{\top}\mB)^{-1}\mB^{\top}\mP\mB\m\Pi. \label{proof1_2}
	\end{align}
	By using the definition of $\m \Omega(\mP) $ given in \eqref{eq:Big-M-SSAProblem-new1}, allows \eqref{eq:OFS-SSAProblem-3} to be replaced by $\m \Pi\m M = \m \Omega(\mP)\m\Pi$. Next, consider $\m \Pi_i\m M_{ij} = \m \Omega(\mP)_{ij}\m\Pi_j$ as the $(i,j)$ element of $\m \Pi\m M = \m \Omega(\mP)\m\Pi$. Then, we have
	\begin{align}
		\m \Pi_i\m M_{ij} = \m \Omega(\mP)_{ij}\m\Pi_j &\Leftrightarrow \m \Pi_i\m M_{ij} = \m\Pi_j\m \Omega(\mP)_{ij} \nonumber\\
		&\Leftrightarrow \m \Pi_i\m M_{ij}-\m\Pi_j\m \Omega(\mP)_{ij}= 0, \label{proof1b}
	\end{align}
	such that
	\[
	\eqref{proof1b} \Leftrightarrow
	\begin{cases}
	\m M_{ij} = \m \Omega(\mP)_{ij}, & \text{if}\;\;\pi_i=1,\, \pi_j=1\\
	\m M_{ij} = 0, \m \Omega(\mP)_{ij}\in\mathbb{R}, &\text{if}\;\;\pi_i=1 ,\,\pi_j=0\\
	\m M_{ij}\in\mathbb{R}, \m \Omega(\mP)_{ij}=0, &\text{if}\;\;\pi_i=0 ,\,\pi_j=1\\
	\m M_{ij}, \m \Omega(\mP)_{ij}\in\mathbb{R}, &\text{if}\;\;\pi_i=0 ,\,\pi_j=0.
	\end{cases}
	\]
	for all $i,j = 1,\hdots,n_u$. For an appropriate large constant $L_2$, the above is equivalent to
	\begin{subequations}
		\begin{align*}
			\vert\m M_{ij}\vert&\leq L_2(1-\pi_i+\pi_j) \\
			\vert\m \Omega(\mP)_{ij}\vert&\leq L_2(1+\pi_i-\pi_j)\\
			\vert\m M_{ij}-\m \Omega(\mP)_{ij}\vert&\leq L_2(2-\pi_i-\pi_j),
		\end{align*}
	\end{subequations}
	that can be written as $\m\Psi_2(\mM,\m \Omega(\mP)) \leq L_2\m\Delta_2(\m\Pi)$ in which $\m\Psi_2(\mM,\m \Omega(\mP))$ and $\m\Delta_2(\m\Pi)$ are given in \eqref{eq:bigM_aux_proof_1c} and \eqref{eq:bigM_aux_proof_1d}.   
	This establishes \eqref{eq:Big-M-SSAProblem-3-2}. Finally, since the left-hand side of \eqref{eq:OFS-SSAProblem-3} and \eqref{proof1_1a} are equal, then \citep[Theorem 2.3.1]{Skelton2013}
	\begin{align}
		\mP\mB\m\Pi &= \mB\mB^{\dagger}\mP\mB\m\Pi \Leftrightarrow
		\mP\mB\m\Pi = \mB(\mB^{\top}\mB)^{-1}\mB^{\top}\mP\mB\m\Pi \nonumber \\
		\Leftrightarrow \mO &= (\mI-\mB(\mB^{\top}\mB)^{-1}\mB^{\top})\mP\mB\m\Pi. \label{proof1c}
	\end{align}
	By using the definition of  $\m \Xi(\mP) $ as in \eqref{eq:Big-M-SSAProblem-new4}, yields $\m \Xi(\mP)\m\Pi = \mO$. Let $\m \Xi(\mP)_{ij}\m\Pi_j$ be the $(i,j)$ element of $\m \Xi(\mP)\m\Pi $. Then, this constraint is equivalent to 
	\[
	\m{\Xi}(\mP)_{ij}\m\Pi_j = 0 \Leftrightarrow 
	\begin{cases}
	\m{\Xi}(\mP)_{ij} = 0, &\text{if}\;\;\pi_j = 1 \\
	\m{\Xi}(\mP)_{ij} \in \mathbb{R}, &\text{if}\;\;\pi_j = 0,
	\end{cases}
	\]
	for $i = 1,\hdots,n_x$ and $j = 1,\hdots,n_u$. For an appropriate large constant $L_3$, the above is equivalent to $\vert\m \Xi(\mP)_{ij}\vert\leq L_3(1-\pi_j)$ such that we obtain $\m\Psi_3(\m \Xi(\mP)) \leq L_3\m\Delta_3(\m\Pi)$, where $\m\Psi_3(\m \Xi(\mP)) $ and $\m\Delta_3(\m\Pi)$ are given in \eqref{eq:bigM_aux_proof_1e} and \eqref{eq:bigM_aux_proof_1f}. This establishes \eqref{eq:Big-M-SSAProblem-3-3}.
	
	The equivalence between \eqref{eq:OFS-SSAProblem} and \eqref{eq:Big-M-SSAProblem} is now summarized. For any feasible $\left\lbrace\m\pi,\m\gamma,\m N,\m M,\m P\right\rbrace $ that satisfies \eqref{eq:OFS-SSAProblem}, by constructing $\m\Theta$ such that $\m\Psi_1(\mN,\m\Theta) \leq L_1\m\Delta_1(\m \Gamma,\m\Pi)$ for a sufficiently large $L_1$, we get $\m{\Pi N \Gamma} = \m\Theta$. Substituting $\m{\Pi N \Gamma} = \m\Theta$ into \eqref{eq:OFS-SSAProblem-2} yields  \eqref{eq:Big-M-SSAProblem-2}. Next, since $\mB$ is full column rank, $\mB^{\top}\mB$ is nonsingular. By using the Moore-Penrose pseudoinverse of $\mB$ given as $\mB^{\dagger}=(\mB^{\top}\mB)^{-1}\mB^{\top}$, pre-multiplying both sides of \eqref{eq:OFS-SSAProblem-3} with $\mB\mB^{\dagger}$ yields \eqref{proof1_2}. By computing $\m \Omega(\mP)$ using \eqref{eq:Big-M-SSAProblem-new1} such that $\m \Pi\m M=\m \Omega(\mP)\m\Pi$, we get $\mM$ and $\m \Omega(\mP)$ satisfy $\m\Psi_2(\mM,\m \Omega(\mP)) \leq L_2\m\Delta_2(\m\Pi)$ for a sufficiently large $L_2$. Then, $\m \Xi(\mP)$ can be computed as \eqref{eq:Big-M-SSAProblem-new4}. Since we have \eqref{proof1c}, $\m \Xi(\mP)$ must satisfy $\m\Psi_3(\m \Xi) \leq L_3\m\Delta_3(\m\Pi)$ for a large constant $L_3$. Therefore, $\left\lbrace\m\pi,\m\gamma,\m N,\m M,\m P,\m\Theta,\m \Omega(\mP),\m \Xi(\mP)\right\rbrace $ is feasible for \eqref{eq:Big-M-SSAProblem}. 
	Conversely, given sufficiently large constants $L_1$, $L_2$, and $L_3$, we always have $\m\Theta = \m{\Pi N \Gamma}$, $\m \Pi\m M=\m \Omega(\mP)\m\Pi$ with $\m \Omega(\mP)$ satisfying \eqref{eq:Big-M-SSAProblem-new1}, and $\m \Xi(\mP)\m\Pi = \mO$ with $\m \Xi(\mP)$ satisfying \eqref{eq:Big-M-SSAProblem-new4} for any feasible $\left\lbrace\m\pi,\m\gamma,\m N,\m M,\m P,\m\Theta,\m \Omega(\mP),\m \Xi(\mP)\right\rbrace $ that satisfies \eqref{eq:Big-M-SSAProblem}. Substituting $\m\Theta = \m{\Pi N \Gamma}$ into  \eqref{eq:Big-M-SSAProblem-2} and \eqref{eq:Big-M-SSAProblem-new1} into $\m \Pi\m M=\m \Omega(\mP)\m\Pi$ yield \eqref{eq:OFS-SSAProblem-2} and $\m\Pi\mM = (\mB^{\top}\mB)^{-1}\mB^{\top}\mP\mB\m\Pi$. The fact that $\mB$ being full column rank implies $(\mB^{\top}\mB)^{-1}\mB^{\top}=\mB^{\dagger}$ so that we have $\m\Pi\mM = \mB^{\dagger}\mP\mB\m\Pi$. Then, pre-multiplying both sides of $\m\Pi\mM = \mB^{\dagger}\mP\mB\m\Pi$ with $\mB$ yields $\mB\m\Pi\mM = \mB\mB^{\dagger}\mP\mB\m\Pi$. Then, substituting \eqref{eq:Big-M-SSAProblem-new4} into $\m \Xi(\mP)\m\Pi=\mO$ yields \eqref{proof1c}, which implies that $\mB\mB^{\dagger}\mP\mB\m\Pi = \mP\mB\m\Pi$ and finally $\mB\m\Pi\mM = \mP\mB\m\Pi$, which is \eqref{eq:OFS-SSAProblem-3}. Hence, $\left\lbrace\m\pi,\m\gamma,\m N,\m M,\m P\right\rbrace $ is feasible for \eqref{eq:OFS-SSAProblem}. This completes the proof.
\end{proof}
\normalcolor
\section{Proof of Theorem \ref{thm:thm2}}\label{apx:D}
\color{black}
\begin{proof}	
	Assume that $\mathbfcal{V}\neq\emptyset$ (that is, SSASP has a solution) and $\mathbfcal{S}_{p}$ be the candidate set at iteration $p$ and $\mathcal{S}_q \in \mathbfcal{S}_p$ with $q = \lceil \sigma/2 \rceil$ and $\sigma=\vert\mathbfcal{S}_p\vert$. Also, let $\mathcal{S}^*_p$ be the best known solution at iteration $p$. If $\mathcal{S}_q$ is infeasible for \eqref{eq:OFS}, then $\mathcal{S}^*_p=\mathcal{S}^*_{p-1}$ and, considering that for practical systems, if a set of selected sensors or actuators renders \eqref{eq:OFS} infeasible, then a subset thereof should also render \eqref{eq:OFS} infeasible, the candidate set is updated such that $\mathbfcal{S}_{p+1}=\mathbfcal{S}_p\setminus\lbrace\mathcal{S}\in\mathbfcal{S}_p\,\vert\,\mathcal{S}_q\vee\mathcal{S} = \mathcal{S}_q \rbrace$. However, assume that $\mathcal{S}_q$ is feasible for \eqref{eq:OFS}, then according to Algorithm \ref{alg:DaC}, $\mathcal{S}^*_p=\mathcal{S}_q$. In this case, for all $\mathcal{S}\in\mathbfcal{U}_p$ where $\mathbfcal{U}_p\triangleq\lbrace\mathcal{S}\in\mathbfcal{S}_p\,\vert\,\mathcal{H}(\mathcal{S}) \geq \mathcal{H}(\mathcal{S}_q)\rbrace$, we have $\mathcal{H}(\mathcal{S}^*_p) \leq \mathcal{H}(\mathcal{S})$. However, since $\mathbfcal{V}_p\subseteq\mathbfcal{U}_p$ where $\mathbfcal{V}_p\triangleq\lbrace \mathcal{S}\in\mathbfcal{U}_p \,|\, \mathcal{S} \textit{ is feasible for \eqref{eq:OFS}}\rbrace$, we have $\mathcal{H}(\mathcal{S}^*_p) \leq \mathcal{H}(\mathcal{S})$ for all $\mathcal{S}\in\mathbfcal{V}_p$. Then, the candidate set is updated such that $\mathbfcal{S}_{p+1} = \mathbfcal{S}_p\setminus\mathbfcal{U}_p$ and the algorithm proceeds. Accordingly, $\sigma$ and $q$ are updated such that $\sigma=\vert\mathbfcal{S}_{p+1}\vert$ and $q = \lceil \sigma/2 \rceil$. If $\mathcal{S}_q$ is infeasible for \eqref{eq:OFS}, where $\mathcal{S}_q\in\mathbfcal{S}_{p+1}$, then $\mathcal{S}^*_{p+1}=\mathcal{S}^*_{p}$. Nonetheless, if $\mathcal{S}_q$ is feasible for \eqref{eq:OFS}, then according to Algorithm \ref{alg:DaC}, $\mathcal{S}^*_{p+1}=\mathcal{S}_q$. 
	
	In the latter case, we have the fact that $\mathcal{H}(\mathcal{S}^*_{p+1})< \mathcal{H}(\mathcal{S}^*_p)$ since $\mathcal{H}(\mathcal{S})<\mathcal{H}(\hat{\mathcal{S}})$ for all $\mathcal{S}\in\mathbfcal{S}_{p+1}$ and $\hat{\mathcal{S}}\in\mathbfcal{S}_{p}$. This shows that at any iteration we have $\mathcal{H}(\mathcal{S}^*_{p+1})\leq \mathcal{H}(\mathcal{S}^*_p)$. Let $l$ denote the index of last iteration. This implies that the sequence $(\mathcal{H}(\mathcal{S}^*_{p}))^l_{p=1}$ is decreasing. Therefore, when $\mathbfcal{S}_l = \emptyset$, $\mathcal{H}(\mathcal{S}^*_l)\leq \mathcal{H}(\mathcal{S})$ for all $\mathcal{S}\in\mathbfcal{V}$. This shows that Algorithm \ref{alg:DaC} computes an optimal solution of SSASP.
\end{proof}
\normalcolor


\end{document}

%% file: preamble.tex
\usepackage{epsfig,color,amsmath,cite}
\usepackage{amsthm}
\usepackage{amsmath}    
\usepackage{bm}
\usepackage{epstopdf}
\usepackage{threeparttable}
\usepackage{amssymb}
\usepackage{url}
\usepackage{enumitem}
\usepackage{multirow}
\usepackage{hhline}
\usepackage{booktabs}
\usepackage[linesnumbered,lined,boxed,commentsnumbered,ruled,longend]{algorithm2e}
\usepackage{comment}
\setlist[itemize]{leftmargin=*}

\DeclareMathOperator*{\minimize}{minimize}

\DeclareMathOperator*{\subjectto}{subject\ to}

\makeatother
\DeclareMathAlphabet\mathbfcal{OMS}{cmsy}{b}{n}

\newtheorem{theorem}{Theorem}
\newtheorem{mydef}{Definition}
\newtheorem{mylem}{Lemma}

\newtheorem{rem}{Remark}
\newtheorem{asmp}{Assumption}

\newtheorem{myprs}{Proposition}
\newtheorem{exmpl}{Example}



\usepackage{stackengine}
\newcommand\barbelow[1]{\stackunder[1.2pt]{$#1$}{\rule{.8ex}{.075ex}}}

\newcommand{\mat}[1]{\boldsymbol{#1}}

\newcommand{\bmat}[1]{\begin{bmatrix} #1 \end{bmatrix}}

\providecommand{\mA}{\ensuremath{\mat{A}}}
\providecommand{\mB}{\ensuremath{\mat{B}}}
\providecommand{\mC}{\ensuremath{\mat{C}}}

\providecommand{\mF}{\ensuremath{\mat{F}}}

\providecommand{\mI}{\ensuremath{\mat{I}}}

\providecommand{\mM}{\ensuremath{\mat{M}}}
\providecommand{\mN}{\ensuremath{\mat{N}}}
\providecommand{\mO}{\ensuremath{\mat{O}}}
\providecommand{\mP}{\ensuremath{\mat{P}}}

\providecommand{\mX}{\ensuremath{\mat{X}}}
\providecommand{\mY}{\ensuremath{\mat{Y}}}

\newcommand{\st}{{\rm s.t.}}

\newcommand{\m}{\boldsymbol}
\allowdisplaybreaks[4]
\usepackage[colorlinks = true,
linkcolor = blue,
urlcolor  = blue,
citecolor = blue,
anchorcolor = blue]{hyperref}

\newcommand{\mc}[1]{\mathcal{#1}}

\usepackage[framemethod=TikZ]{mdframed}
\mdfdefinestyle{MyFrame}{%
	linecolor=black,
	outerlinewidth=1pt,
	roundcorner=1pt,
	innerrightmargin=5pt,
	innerleftmargin=5pt,}
	
\usepackage{graphicx}

\usepackage{cleveref}